\documentclass[11pt,reqno]{amsart}
\usepackage[T1]{fontenc}
\usepackage{amsmath,amssymb,amsthm,lmodern}
\usepackage{mathtools}
\usepackage[margin=1in]{geometry}
\usepackage{microtype,needspace}
\usepackage{tikz-cd}
\usetikzlibrary{positioning,arrows.meta,calc}
\usepackage{hyperref}
\hypersetup{hidelinks,pdftitle={Canonical Projectivization of Algebraic Varieties},pdfauthor={Parsa Bakhtary}}
\theoremstyle{plain}
\newtheorem{theorem}{Theorem}[section]
\newtheorem{proposition}[theorem]{Proposition}
\newtheorem{lemma}[theorem]{Lemma}
\newtheorem{corollary}[theorem]{Corollary}
\theoremstyle{definition}
\newtheorem{definition}[theorem]{Definition}
\newtheorem{example}[theorem]{Example}
\theoremstyle{remark}
\newtheorem{remark}[theorem]{Remark}
\newcommand{\OO}{\mathcal O}
\newcommand{\PP}{\mathbb P}
\newcommand{\MM}{\mathcal M}
\newcommand{\QQ}{\mathbb Q}

\newcommand{\m}{\mathfrak m}
\newcommand{\Bl}{\operatorname{Bl}}
\newcommand{\nBl}{\operatorname{Bl}^{\nu}}
\DeclareMathOperator{\Spec}{Spec}
\DeclareMathOperator{\Supp}{Supp}
\DeclareMathOperator{\length}{length}
\DeclareMathOperator{\Aut}{Aut}
\DeclareMathOperator{\Hilb}{Hilb}
\begin{document}
\title[Canonical projectivization]{Canonical Projectivization of Algebraic Varieties}
\author[P.~Bakhtary]{Parsa Bakhtary}
\address{Independent, Sunnyvale, California, USA}
\email{pbakhtary@gmail.com}
\date{}

\begin{abstract}
We construct canonical projective birational models of normal proper
integral algebraic spaces over arbitrary fields. The construction preserves
the intrinsic quasi-projective core and is natural under semilinear
isomorphisms. For smooth proper integral input in characteristic zero, we
combine an intrinsic choice of ideals with W{\l}odarczyk's principalization
to obtain a canonical tower of blowups in smooth centers ending in a smooth
projective variety. Consequently, every compact connected Moishezon manifold
admits a smooth projective bimeromorphic model to which all automorphisms
lift compatibly. The main ingredient is a family-sum argument yielding an
intrinsic finite selection of maximal tuples of boundary ideal germs;
successive modifications strictly decrease the boundary-image dimension. We
also construct canonical quasi-projective models for nonproper normal input
and a greatest schematizing ideal for separated finite-type algebraic spaces
over noetherian bases. Oda's threefold illustrates the construction.
\end{abstract}

\maketitle

\section{Introduction}\label{sec:introduction}

Chow's lemma provides projective modifications of proper algebraic
varieties, and its open-preserving form can leave a prescribed
quasi-projective open unchanged. These constructions involve choices of
opens, compactifications and modifying ideals. We construct the ideals
intrinsically, obtaining a canonical projective model over every field.
In characteristic zero, smooth input also admits a canonical tower of
smooth blowups with a smooth projective model. An intrinsic construction
lets every automorphism of the original space lift to the same model;
the resulting group action is one reason to seek canonicity.

For a normal integral separated algebraic space $X$ of finite type over
a field $k$, let $\mathcal Q(X)$ be its set of maximal
$k$-quasi-projective open subsets. Its \emph{quasi-projective core} is
\[
 U(X)=\bigcap_{W\in\mathcal Q(X)}W.
\]
The finiteness theorem of W{\l}odarczyk
\cite[Theorem~A]{Wlodarczyk1999}, in Benoist's generality
\cite[Theorem~9]{Benoist}, makes this a nonempty quasi-projective open.
For algebraic spaces, apply the theorem to the largest schematic open;
see Section~\ref{sec:prelim}. This is the original open preserved by
our construction.

A semilinear isomorphism consists of an isomorphism of ground fields
and an isomorphism of the given spaces after transport of scalars.
A normalized blowup means the normalization of an ordinary blowup.
In the normal construction, one round is one normalized blowup.
In the smooth construction, a round is one application of the fixed
principalizer and may contain several smooth blowups.

\begin{theorem}[Normal projective models]\label{thm:normal}
For every field $k$ and every normal integral separated algebraic space
$X$ of finite type over $k$, there is a canonical finite ordered tower
of normalized blowups
\[
 P(X)=X_r\longrightarrow X_{r-1}\longrightarrow\cdots\longrightarrow X_0=X
\]
with the following properties.
\begin{enumerate}
\item The model $P(X)$ is a normal quasi-projective variety. If $X$ is
proper, $P(X)$ is projective over $k$.
\item Every map is an isomorphism over the inverse image of the original
core $U(X)$.
\item The tower is empty if $X$ is quasi-projective. It has at most
$\dim X$ nonidentity rounds for scheme input, and at most $\dim X+1$
for algebraic-space input.
\item Semilinear isomorphisms induce unique compatible isomorphisms of
the entire ordered towers, respecting identities and composition.
\end{enumerate}
\end{theorem}

No perfectness, geometric-normality or geometric-integrality hypothesis
is imposed. The morphism $P(X)\to X$ is projective and birational even
when $X$ is not proper. In that case we call $P(X)$ the canonical
quasi-projective model. The marked scheme version,
Theorem~\ref{thm:marked}, preserves any specified nonempty
quasi-projective open.

\begin{theorem}[Smooth projective models]\label{thm:smooth}
Fix the principalization procedure and divisor convention in
Section~\ref{sec:principalizer}. For every characteristic-zero field $k$
and every smooth proper integral algebraic space $X$ over $k$, there is
a canonical finite ordered tower of blowups in smooth closed subspaces
\[
 P_{\mathrm{sm}}(X)=X_r\longrightarrow\cdots\longrightarrow X_0=X
\]
with smooth projective variety $P_{\mathrm{sm}}(X)$. All centers avoid
the inverse image of $U(X)$, the tower is empty on projective input,
and semilinear isomorphisms lift uniquely through the entire tower,
compatibly with identities and composition.
\end{theorem}

These are two constructions with a common selection argument. Normalized
blowups give the arbitrary-field result; principalization supplies smooth
centers in characteristic zero. Canonicity means naturality under the
stated isomorphisms. The ideals are defined by intrinsic conditions
on the given space and its retained tower; no representatives of
isomorphism classes enter their definition.
Neither construction asserts compatibility with
recomputation after arbitrary field extension or under smooth morphisms.
Both are non-effective: existence conditions and intrinsic minima define
the selected ideals, without an algorithm for deciding those conditions.

\begin{corollary}\label{cor:automorphisms}
For a normal proper integral algebraic space $X$ over $k$, lifting gives
an injective homomorphism of abstract groups
\[
 \Aut_k(X)\lhook\joinrel\longrightarrow\Aut_k(P(X)).
\]
For smooth proper $X$ in characteristic zero, the same holds with the
smooth projective model $P_{\mathrm{sm}}(X)$. The models have the same
dimension as $X$.
\end{corollary}

\begin{proof}
The lifts respect composition. A lift equal to the identity forces the
original automorphism to be the identity on the dense open $U(X)$,
and hence on the integral separated space $X$.
\end{proof}

Here lifting concerns the abstract group of automorphisms: each
automorphism has a unique lift, and the lifts respect composition.
No compatibility with families of automorphisms is asserted. The
algebraic torus actions used in Section~\ref{sec:oda} will follow from
an additional invariance argument.

The passage from normal algebraic spaces to schemes uses a companion
construction over a more general base. Call a coherent ideal on an
algebraic space \emph{schematizing} if it is unit on the largest
schematic open and its ordinary blowup is a scheme.

\begin{theorem}[Canonical schematization]\label{thm:schematization}
Let $S$ be a noetherian scheme and $X\to S$ a separated morphism of
finite type from an algebraic space, possibly nonreduced. There is a
greatest schematizing ideal $H_X\subseteq\OO_X$. Its blowup is a
scheme, projective over $X$, and unchanged over the entire schematic
locus $X_{\mathrm{sch}}$. The ideal and its blowup are natural under
isomorphisms of $X$, independently of the chosen presentation over $S$.
The construction is the identity when $X$ is a scheme.
\end{theorem}

Greatestness refers to inclusion of ideals. Theorem~\ref{thm:schematization}
does not give a least model or require a field, reducedness or normality.
Its proof in Section~\ref{sec:schematization} is independent of the
Hilbert-scheme selection argument.

\subsection{Earlier work and the equivariant question}

Rydh explains in his discussion of Chow's lemma \cite[\S3.3]{Rydh}
that the choice of a quasi-projective open and a projective
compactification prevents that construction from being functorial.
He asks to what extent a functorial Chow lemma exists. Theorems~\ref{thm:normal}
and~\ref{thm:smooth} give isomorphism-natural constructions while
preserving the intrinsic open $U(X)$.

The natural Chow lemma for smooth threefolds and fourfolds
\cite[Theorem~1.1]{Bakhtary} uses numerical data on maximal-open
boundaries. Its final step treats boundary images of dimension at most
one. Here ideal germs at the generic points of those images replace
that numerical argument and give a construction in every dimension.
For toric varieties, the earlier rule \cite[\S5]{BakhtaryToric}
uses wall characters and is natural under lattice isomorphisms.
Section~\ref{sec:oda} compares the two constructions on Oda's threefold.

Meng, Perroni and Zhang \cite[Question~1.4]{MPZ} ask whether every
compact Moishezon manifold admits a projective bimeromorphic model to
which its full automorphism group lifts, and ask the corresponding
K\"ahler-model question for Fujiki's class $\mathcal C$.
Corollary~\ref{cor:moishezon} answers the Moishezon question for the
abstract automorphism group, with a smooth blowup tower preserving the
quasi-projective core. Jia and Meng
\cite[Theorem~1.1 and Remark~1.2]{JiaMeng} construct a K\"ahler model
equivariant under the subgroup preserving a specified big $(1,1)$-class.
Their proof sums translates of an ideal and uses residual ideals whose
sum is the unit ideal \cite[\S3, proof of Theorem~1.1 and Claim~3.3]{JiaMeng}.
Our Moishezon result requires no invariant big class and no finiteness
assumption on the action on cohomology. Example~\ref{ex:infinite-cohomology}
gives a smooth proper nonprojective fivefold with an automorphism of
infinite order on $H^2$ and with no invariant big class, even for the
cyclic subgroup it generates.

Kim \cite[Theorem~1.1]{Kim} states the full-automorphism K\"ahler-model
conclusion for Fujiki's class $\mathcal C$. The proof on p.~564 contains
an incorrect equivariant-flattening assertion: after choosing a K\"ahler
surjection, it asserts that flattening can be achieved while every
automorphism of the target lifts. This assertion fails already when the chosen map is the
blowup of a point of $\PP^2$; Remark~\ref{rem:kim-gap} gives the
counterexample. As written, the published argument does not establish the
stated theorem; the counterexample does not disprove its existential
conclusion. Our proof establishes the full-automorphism conclusion
in the Moishezon case independently, together with canonicity and
preservation of the original core. We do not address the general
class-$\mathcal C$ question.

\subsection{The selection mechanism}

Let $f\colon Y\to X$ be the modification built so far, and let $W$ be
a maximal quasi-projective open containing the preserved mark. To improve
the model, we lower the dimension of the closed image $f(Y\setminus W)$.
It suffices to find, after modification, one quasi-projective open
containing the lifted $W$ and the whole fibers over the generic points
of the largest-dimensional components of that image. We call tuples of
ideal germs \emph{enveloping} when their contracted ideal has a normalized
blowup with such a neighborhood.

The existence of a witness follows from the open-preserving Chow lemma.
To make an intrinsic selection, we first take maximal enveloping tuples
under ideal inclusion. We then minimize their thickness and a bound on
the coefficients of their length-Hilbert polynomials. The resulting
collection lies in a noetherian parameter space, but boundedness alone
does not make it finite. The additional fact is that the sum of a dense
irreducible family of enveloping tuples is enveloping again. Full
maximality forces every member of such a family to equal its sum, leaving
only finitely many selected tuples.

The family-sum argument is the main point of the proof. On the blowup
of the sum, failure of one family member to generate the invertible
image ideal at a fixed fiber point is a proper closed condition on
parameters. Irreducibility lets one avoid finitely many such conditions
simultaneously. A member avoiding them transfers its quasi-projective
neighborhood to the chosen finite set. Finiteness of maximal
quasi-projective opens then turns this finite-point statement into one
neighborhood of whole fibers. We spell out these two distinct finiteness
steps in Section~\ref{sec:families}.

Take the product of all selected ideals. Making that product invertible
makes every factor invertible, so the resulting modification maps to
each selected witness. Their neighborhoods then give a strict drop in
the boundary-image dimension. Normalized blowup or smooth principalization
performs this step; the drop proves termination in both cases.

Sections~\ref{sec:prelim}--\ref{sec:construction} contain the scheme
construction: local ideal data, finite selection, and dimension drop.
For a first reading, the definition of an enveloping tuple in
Section~\ref{sec:germs}, the finite-selection argument in
Section~\ref{sec:selection}, and the dimension drop in
Section~\ref{sec:construction} give the structure of the proof.
Section~\ref{sec:oda-core} computes these objects on Oda's threefold,
where each maximal open has a unique maximal enveloping germ.
The detailed family argument in Sections~\ref{sec:bad-parameters}
and~\ref{sec:family-sum} explains why selection is finite.

Section~\ref{sec:schematization} constructs a scheme model and applies
the marked scheme theorem to normal algebraic spaces.
Section~\ref{sec:smooth} gives the smooth construction, first for
schemes and then directly for algebraic spaces, and derives the
Moishezon consequence. It uses the same selection at every stage
and does not depend on the schematizing-ideal theorem.
Section~\ref{sec:oda} concludes with the Oda comparison and the example
with infinite cohomological action.

\subsection{Conventions}

A variety over a field is an integral separated scheme of finite type.
Algebraic spaces are separated and of finite type over the indicated
base. Fields are arbitrary unless characteristic zero is explicitly
assumed. Quasi-projectivity and projectivity refer to the ground field
unless a relative base is specified. For an ideal $I$ on a normal variety
$Y$, write $\nBl_IY$ for the normalization of its ordinary blowup.
Normalizations of varieties over fields are finite. A dominant morphism
from a normal variety to a blowup factors uniquely through its
normalization; we use this in Propositions~\ref{prop:pooling}
and~\ref{prop:drop}. Ideal pullbacks
always denote their images in the structure sheaf, rather than tensor
pullbacks. An ideal is \emph{unit} on an open if that image is the
structure sheaf there. Closed points of a generic fiber mean points
closed in that fiber, not necessarily in the total space. We use
$\dim\varnothing=-1$.

\section{Maximal opens and projective modifications}\label{sec:prelim}

W{\l}odarczyk proved finiteness of maximal quasi-projective opens for
complete normal varieties over an algebraically closed field when
$(\operatorname{Div}/\operatorname{Car})\otimes\QQ$ is finite-dimensional,
including smooth varieties in characteristic zero
\cite[Theorem~A and the following remark, p.~41]{Wlodarczyk1999}.
Benoist proves that a separated normal scheme of finite type over
an arbitrary field has finitely many maximal quasi-projective open
subschemes \cite[Theorem~9]{Benoist}. Neither properness nor geometric
normality is required. We need this generality for normalized blowups,
which need not be locally factorial.

Noetherianness places every quasi-projective open in a maximal one.
Their finite intersection $U(X)$ is nonempty by integrality and is
quasi-projective. It equals $X$ precisely when $X$ is quasi-projective. For proper $X$,
this is equivalent to projectivity. The same finiteness holds for a normal integral algebraic space $X$ of
finite type over $k$. Indeed, its largest schematic open $X_{\mathrm{sch}}$
is a normal variety \cite[Tags~03JH and~06NH]{Stacks}. Every quasi-projective open
of $X$ is a scheme and lies in $X_{\mathrm{sch}}$, so the two spaces
have exactly the same maximal quasi-projective opens. In particular,
the core of such an $X$ is again a nonempty quasi-projective scheme.

Finiteness also lets us pass from neighborhoods of
finitely many fiber points to a neighborhood of whole fibers.

\begin{lemma}\label{lem:finite-neighborhoods}
Let $T$ be a normal integral separated algebraic space of finite type
over $k$, let $W\subset T$ be a
quasi-projective open, and let $T\to X$ be a projective morphism. Fix finitely many points $\eta_i\in X$.
Suppose that every finite set of closed points of the fibers $T_{\eta_i}$
is contained, together with $W$, in a quasi-projective open of $T$.
Then some quasi-projective open contains $W$ and every $T_{\eta_i}$.
\end{lemma}

\begin{proof}
List the finitely many maximal quasi-projective opens of $T$ containing
$W$. If none contains all the indicated fibers, choose for each open a
missing closed point of one fiber. Such a point exists because the
intersection of its closed complement with that fiber is a nonempty closed
subset of a fiber of finite type over a field. A quasi-projective open
containing $W$ and all the chosen points is contained in a maximal one,
contradicting the choice of its missing point.
\end{proof}

\subsection{A modification preserving a given open}

\begin{lemma}[Open-preserving Chow lemma, ideal form]\label{lem:chow-ideal}
Let $Y$ be a normal variety over a field $k$, and let $W\subset Y$
be a nonempty quasi-projective open. There is a nonzero coherent ideal
$I\subset\OO_Y$, unit on $W$, such that $\Bl_IY$ is quasi-projective
over $k$. Its normalization is quasi-projective as well. If $Y$ is
proper, both modifications are projective over $k$.
\end{lemma}

\begin{proof}
Choose a proper scheme compactification $Y\subset\overline Y$
\cite[Tag~0F41]{Stacks}, replacing it by the closure of $Y$ if necessary.
The open $W$ is also an open of $\overline Y$. Rydh's open-preserving
Chow lemma gives a finite sequence of $W$-admissible blowups
$Z\to\overline Y$ with $Z$ projective over $k$
\cite[Theorem~3.5]{Rydh}. Here admissibility means that each center
is disjoint from the inverse image of $W$.

A finite composite of such blowups is a single $W$-admissible blowup
\cite[Tag~080L]{Stacks}. Restricting this blowup to $Y$ gives
$\Bl_IY$ for a coherent ideal $I$ unit on $W$. Its source is an open
subscheme of the projective scheme $Z$, hence is quasi-projective.
The ideal is nonzero because it is unit on the dense open $W$.
Normalization is finite for varieties over a field, so
$\nBl_IY$ is quasi-projective too. If $Y$ is proper, both sources
are proper over $k$; properness and quasi-projectivity give projectivity.
\end{proof}

This lemma supplies witnesses for the conditions below. The construction
will not retain the compactification or the witnessing ideal. Their
only role is to show that the collection from which we select is nonempty.

\subsection{Retained divisors and a relative polarization}\label{sec:polarization}

The numerical bounds used for finite selection will require an ample
line bundle on each finite thickening of a fiber. We obtain it from the
ordered tower itself, so it introduces no further choice. It will enter
only through the length-Hilbert polynomial in \eqref{eq:polynomial};
maximal enveloping tuples and their least thickness are defined before
that polynomial is used to break ties.

For the normal construction, a retained tower is an ordered sequence
of normalized blowups
\[
 Y=X_r\longrightarrow X_{r-1}\longrightarrow\cdots\longrightarrow X_0=X.
\]
If the $j$th map is the normalized blowup of $I_j$, retain the effective
Cartier divisor defined by the invertible image ideal
$I_j\OO_{X_j}$. Write $E_j$ for its total transform on $Y$.
The full Cartier divisor is part of the data: replacing it by its
reduced support would change the line bundles used below.

The line bundle $\OO_{X_j}(-E_j)$, before total transform, is ample
relative to $X_{j-1}$. Indeed, it is the finite pullback of the
tautological ample line bundle on the ordinary blowup. Composition of
relative ampleness gives positive integers $a_1,\ldots,a_r$ for which
\begin{equation}\label{eq:relative-polarization}
 A_f=\OO_Y\Bigl(-\sum_{j=1}^r a_jE_j\Bigr)
\end{equation}
is ample relative to the composite $f\colon Y\to X$
\cite[Tag~0C4K]{Stacks}.
Choose the lexicographically least such positive tuple. It exists:
minimize the first coordinate among the nonempty set of allowed tuples,
then the second among those remaining, and continue through the finitely
many coordinates. For the identity tower set $A_f=\OO_X$.

Thus the ordered tower supplies its own relative polarization. No ample
bundle on $X$ or choice of a projective model is needed. The property of
being $f$-ample, the retained divisors, and the least tuple are transported
by isomorphisms of ordered towers. The same definition will apply to the
ordered smooth blowup towers in Section~\ref{sec:smooth}, with their
stated divisor convention. Sections~\ref{sec:germs}--\ref{sec:selection}
use only projectivity of $f$, normality of the varieties, the preserved
open, and the resulting line bundle $A_f$.

\section{Ideal germs at the boundary}\label{sec:germs}

The next round should enlarge a quasi-projective open across the
fibers over the largest components of its boundary image. We first
describe this geometric objective, then the local ideal data that encode
it. The completed rings below are used to parameterize ideals on finite
infinitesimal neighborhoods. They will also serve the smooth
algebraic-space construction (Lemma~\ref{lem:completed-point});
Remark~\ref{rem:scheme-comparison} explains why ordinary local rings
give exactly the same selection for varieties.

Fix a nonempty quasi-projective open $V\subset X$ in a normal variety
over an arbitrary field. Let $f\colon Y\to X$ be a retained tower of
normal varieties, isomorphic over $V$, with relative polarization $A_f$.
Suppose that $Y$ is not quasi-projective. Choose a maximal quasi-projective
open $W\subset Y$ containing $f^{-1}V$ and
write
\[
 B=(Y\setminus W)_{\mathrm{red}},\qquad d_W=\dim f(B).
\]
The image $f(B)$ is closed since $f$ is proper. Let
$D_1,\ldots,D_s$ be its irreducible components of dimension $d_W$, with
generic points $\eta_1,\ldots,\eta_s$. Tuples below are indexed by these components; the numbering makes no
choice. The image is a proper closed subset of $X$ because $f$ is
birational and $W$ contains the generic point of $Y$.

Suppose that a further projective modification admits a quasi-projective
open containing the lifted $W$ and every whole fiber over the $\eta_i$.
The image of its complement in $X$ is closed, lies in $f(B)$, and misses
all the top-dimensional generic points of $f(B)$. Its dimension is
therefore smaller. This is the improvement that our ideal tuples will
encode. Containing one point of a fiber would not suffice: the
complementary image misses a base point exactly when the open contains
its entire fiber. We measure the image in the original $X$, since a
blowup may create exceptional divisors even while this image shrinks.

For each $i$, put
\[
 R_i=\widehat{\OO}_{X,\eta_i},\qquad \kappa_i=k(\eta_i),
\]
and let $\m_i$ be the maximal ideal of $R_i$. Each $R_i$ is a complete noetherian equicharacteristic local ring.
Regularity is not required. The natural map $\Spec R_i\to X$
is flat, being completion followed by localization, and affine, since
$X$ is separated. Set
\[
 Y_i=Y\times_X\Spec R_i,\qquad W_i=W\times_Y Y_i.
\]
Then $Y_i$ is projective over $R_i$, as the base change of the projective morphism $f$. For $b\geq1$, write
\[
 \Lambda_{i,b}=\OO_{X,\eta_i}/\mathfrak m_{X,\eta_i}^{b}
         =R_i/\m_i^b.
\]
Here $\mathfrak m_{X,\eta_i}$ is the maximal ideal of $\OO_{X,\eta_i}$.
These Artinian rings record transverse information up to order $b$.
In particular,
\begin{equation}\label{eq:completion-thickening}
 \Spec R_i\times_X\Spec \Lambda_{i,b}=\Spec \Lambda_{i,b}.
\end{equation}
Indeed, $\Spec\OO_{X,\eta_i}\to X$ is a localization and hence a
monomorphism, so this is the identity
$R_i\otimes_{\OO_{X,\eta_i}}\Lambda_{i,b}=\Lambda_{i,b}$. We will use this
finite-level identity to recover an ideal after contraction. Write
\[
 Y_{i,b}=Y\times_X\Spec\Lambda_{i,b}
        =Y_i\times_{R_i}\Spec\Lambda_{i,b}
\]
for the common finite thickening.

An \emph{ideal germ} at $\eta_i$ is a coherent ideal on $Y_i$.
Since $Y_i$ need not be affine, it is not in general an ideal of $R_i$.
Consider tuples $h=(h_1,\ldots,h_s)$ of nonzero ideal germs that are
unit on $W_i$; we allow the unit ideal $\OO_{Y_i}$.
Their cosupports lie in the closed fibers of $Y_i\to\Spec R_i$.
Indeed, a strict generization of $\eta_i$ in $f(B)$ would contradict
the definition of $d_W$, and the inverse image of $\eta_i$ is the
closed point of $\Spec R_i$.

Coherence and quasi-compactness now imply that
$\OO_{Y_i}/h_i$ is killed by some power $\m_i^b$. Thus $h_i$ is
the inverse image of an ordinary coherent ideal on the projective
finite thickening $Y_{i,b}$: all its relevant data have finite thickness.
For a field $\kappa$, for example, $(x,y^2)\subset\kappa[[x,y]]$ contains $(x,y)^2$
and is recovered from its image modulo that square. The same ideal on
the common quotient recovers $(x,y^2)$ in
$\kappa[x,y]_{(x,y)}$. This example illustrates finite-thickness
recovery only. In general $Y_i$ may have a positive-dimensional fiber, and
the data are ideal sheaves on $Y_i$, not ideals of the base ring $R_i$.
No algebraization of infinite formal data is involved.

To form a blowup on $Y$, we need to turn the tuple of germs into a
global ideal. The useful operation is contraction: retain exactly those
local functions whose images belong to every prescribed germ. With
$j_i\colon Y_i\to Y$, define the \emph{canonical contraction} of $h$ by
\begin{equation}\label{eq:contraction}
 H(h)=\bigcap_{i=1}^s\ker\bigl(\OO_Y\longrightarrow
       (j_i)_*(\OO_{Y_i}/h_i)\bigr).
\end{equation}

The contraction must have two properties: it must recover every chosen
germ, and it must introduce no change on $W$. The first property is what
allows all subsequent tests to be made on the finite neighborhoods.

\begin{lemma}\label{lem:contraction}
The sheaf $H(h)$ is a nonzero coherent ideal, unit on $W$, and
$H(h)\OO_{Y_i}=h_i$ for every $i$. If two coherent ideals on $Y$ are
unit on $W$ and have the same germs on all the $Y_i$, they agree on an
open neighborhood of $W$ and all the fibers $Y_{\eta_i}$.
\end{lemma}

\begin{proof}
The morphisms $j_i$ are affine and flat, as observed above. Their pushforwards are quasi-coherent, so
\eqref{eq:contraction} is a quasi-coherent ideal and is coherent by
noetherianness. Restriction to $W$ makes every displayed quotient zero.

For each $i$ there is an integer $b_i\geq1$ such that
$\m_i^{b_i}\OO_{Y_i}\subset h_i$. This follows from the support
condition and quasi-compactness. If $\mathcal I_{D_i}$ is the ideal of
$D_i$ in $X$, then
\[
 \prod_i(\mathcal I_{D_i}\OO_Y)^{b_i}\subset H(h).
\]
The left side is nonzero since $Y$ is integral and $f$ is dominant.
For $i\ne j$, the $j$th quotient is annihilated by
$(\mathcal I_{D_j}\OO_Y)^{b_j}$. This ideal becomes the unit ideal
at $\eta_i$, since $\eta_i\notin D_j$, so that quotient vanishes
after completed pullback.

Fix $i$, put $b=b_i$, and write $\OO_{Y_i}/h_i=(\iota_i)_*Q_i$,
where $\iota_i\colon Y_{i,b}\hookrightarrow Y_i$ is the closed immersion
and $Q_i$ is coherent on $Y_{i,b}$. Then
$(j_i)_*(\OO_{Y_i}/h_i)=(j_i\iota_i)_*Q_i$. The morphism
$j_i\iota_i\colon Y_{i,b}\to Y$ is the base change of the affine
morphism $\Spec\Lambda_{i,b}\to X$. By
\eqref{eq:completion-thickening},
\[
 Y_i\times_Y Y_{i,b}
 =Y\times_X(\Spec R_i\times_X\Spec\Lambda_{i,b})=Y_{i,b},
\]
with projections $\iota_i$ and the identity. Affine base change gives
$j_i^*(j_i\iota_i)_*Q_i\cong(\iota_i)_*Q_i$. Under this identification,
the counit
$j_i^*(j_i)_*(\OO_{Y_i}/h_i)\to\OO_{Y_i}/h_i$
is the identity. Thus the other components have disappeared, while
the $i$th quotient has returned unchanged. Pulling back the kernel sequence by
the flat map $j_i$ gives exactly $h_i$. Flatness also preserves the
finite intersection in \eqref{eq:contraction}, proving simultaneous
recovery.

Finally, for two ideals $I,J$, the coherent sheaves $(I+J)/I$ and
$(I+J)/J$ measure their difference. Their completed pullbacks vanish.
For each $y\in Y_{\eta_i}$, the common closed fiber gives a point
$y'\in Y_i$ above $y$. The local map
$\OO_{Y,y}\to\OO_{Y_i,y'}$ is flat and local, hence faithfully flat.
The two difference quotients therefore have zero stalk at $y$.
Their supports miss every $Y_{\eta_i}$ and also $W$. The complement
of those supports is the required neighborhood.
\end{proof}

\begin{definition}\label{def:enveloping}
A tuple $h$ is \emph{enveloping} if the normalized blowup
\[
 T_h=\nBl_{H(h)}Y
\]
has a $k$-quasi-projective open containing the inverse image of $W$ and
all points of every fiber $(T_h)_{\eta_i}$.
Maximality of an enveloping tuple will always mean maximality under componentwise ideal inclusion among \emph{all}
enveloping tuples for the fixed $f,W$.
\end{definition}

An enveloping tuple need not make all of $T_h$ quasi-projective. It
requires a single quasi-projective neighborhood of the old open and the
specified whole fibers. This is exactly the improvement needed to remove
the largest components from the next boundary image. In
Section~\ref{sec:oda-core}, the identity tower on Oda's threefold gives
the simplest case: $Y_i=\Spec R_i$, and the selected germ is the
maximal ideal at the generic point of one missing curve.

The neighborhood in this definition is quasi-projective over $k$;
projectivity over the local bases $\Spec R_i$ alone would not
supply it. Nevertheless, the condition depends only on the indicated
germs. Indeed, by Lemma~\ref{lem:contraction}, any other global
ideal which is unit on $W$ and has these germs agrees with $H(h)$
near $W$ and all the fibers. The normalized blowups agree there
as well, so we can restrict a witnessing open to that neighborhood
and transport it from one blowup to the other.

\begin{lemma}\label{lem:existence}
Enveloping tuples exist, and every enveloping tuple is contained in a maximal enveloping
tuple.
\end{lemma}

\begin{proof}
Choose an ideal $I$ with quasi-projective blowup, unit on $W$, by the
open-preserving Chow lemma (Lemma~\ref{lem:chow-ideal}). Its completed pullbacks have the required support. The quasi-projective
normalized blowup of $I$ gives the required
neighborhood, which transports to the contracted ideal by the preceding
paragraph. Ascending chains of tuples of coherent ideals stabilize on
the noetherian quasi-compact spaces $Y_i$. The
same ascending-chain condition holds in the enveloping subposet and gives
maximal elements above every member.
\end{proof}

\section{Families of ideal germs}\label{sec:families}

We would like to turn a bounded collection of enveloping tuples
into a finite collection. Hilbert schemes allow us to place the
ideals in parameter spaces of finite type, but this still leaves
open the possibility of infinitely many rational points. The
additional fact we need is that the sum of a dense irreducible
family of enveloping tuples is again enveloping. Once this is
proved, maximality will force each member of such a family to
equal the sum.

\begin{figure}[tb]
\centering
\begin{tikzpicture}[>=stealth, node distance=8mm,
  box/.style={draw=black!45,rounded corners=3pt,fill=blue!3,
    text width=0.86\linewidth,align=center,inner sep=7pt,font=\small}]
\node[box] (finite) {Fix a finite set $\Sigma$ of fiber points.\\
Each $\tau\in\Sigma$ excludes a proper closed subset $F_\tau\subset C$.};
\node[box,below=of finite] (choose) {Irreducibility and density give
$a\in\mathcal A\setminus\bigcup_{\tau\in\Sigma}F_\tau$.\\
Its neighborhood transfers to a quasi-projective open containing
$W\cup\Sigma$.};
\node[box,below=of choose,fill=green!3] (whole) {There are only finitely many
maximal quasi-projective opens containing $W$.\\
If each missed a fiber point, those finitely many points would contradict
the preceding conclusion.};
\draw[->,black!65] (finite)--(choose);
\draw[->,black!65] (choose)--(whole);
\end{tikzpicture}
\caption{The two finiteness arguments in the family-sum proof. The chosen
parameter may depend on $\Sigma$. A neighborhood of all the fibers is
obtained only in the last step, using finiteness of maximal opens.}
\label{fig:finite-avoidance}
\end{figure}
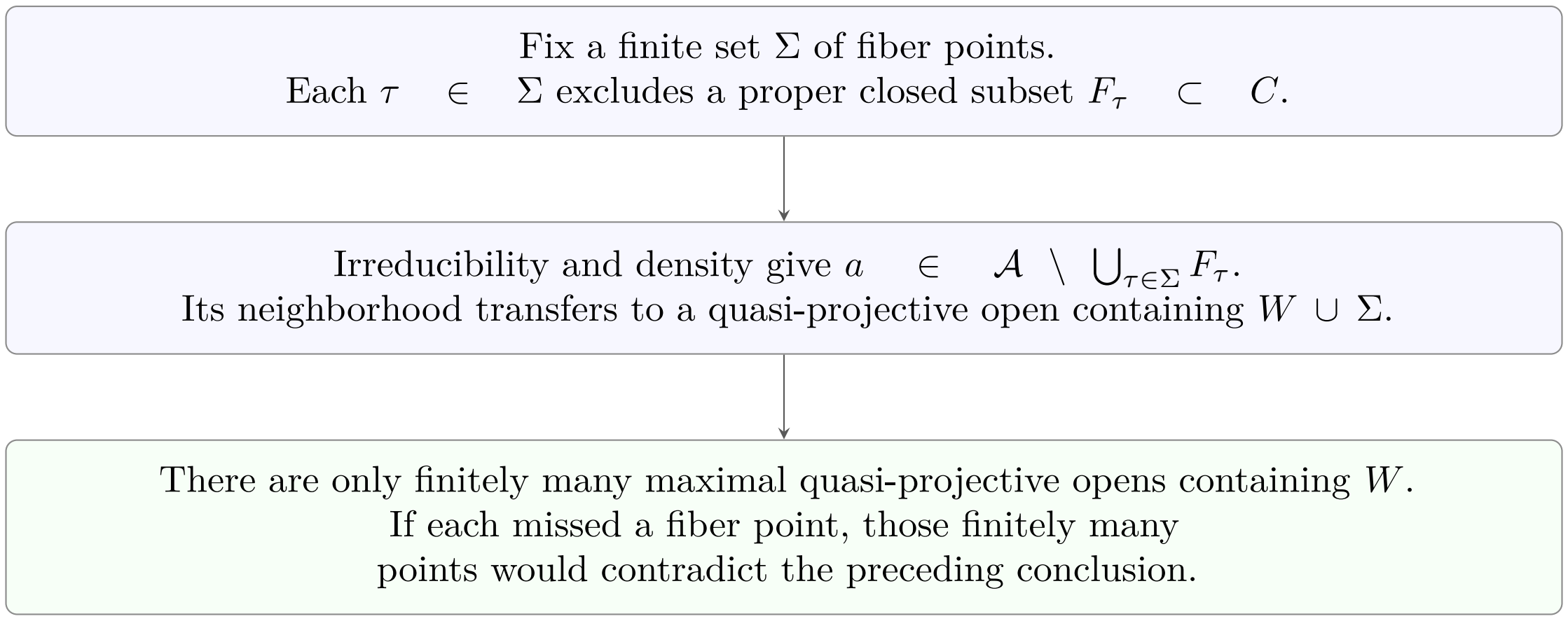

To understand why the sum has the required neighborhood, we look
at its blowup and ask where a member of the family fails to
generate it. At a closed point of a generic fiber, this failure
will be a proper closed condition on the parameters. We can avoid
finitely many such conditions at once, and then use
Lemma~\ref{lem:finite-neighborhoods} to obtain a neighborhood of
the whole fibers. Figure~\ref{fig:finite-avoidance} separates these two steps. We carry
them out after describing the parameter spaces. Throughout this section, the data $f,W,D_i$ are those of
Section~\ref{sec:germs}.

\subsection{Parameter spaces}\label{sec:parameters}

Recall the part of the Hilbert scheme theorem that we use. If $Z$ is
a projective scheme over a field $F$, with a very ample line bundle
$L$, then $\Hilb^P(Z)$ is a projective $F$-scheme whose $F$-points
correspond to closed subschemes $Z'\subset Z$ with Hilbert polynomial
$P(t)=\chi(Z',L^{\otimes t}|_{Z'})$. There is a universal closed
subscheme of $Z\times_F\Hilb^P(Z)$, flat over $\Hilb^P(Z)$.
In this way, a point of the Hilbert scheme records a subscheme
inside the fixed $Z$, including its defining ideal. This is the
information needed when we later take sums of ideals. For a
zero-dimensional $Z'$, its Hilbert polynomial is simply the
constant $\dim_F H^0(Z',\OO_{Z'})$; for higher-dimensional
subschemes it records the eventual growth of sections.

To apply this description to our ideals over local rings, we first
put a bound on their thickness. This places them on a projective
scheme over an Artinian ring. Choosing a coefficient field then
puts us in the situation of the Hilbert scheme theorem.
Fix $b\geq1$. Recall that $Y_{i,b}$ is the closed subscheme of $Y_i$
defined by $\m_i^b\OO_{Y_i}$, projective over $\Lambda_{i,b}=R_i/\m_i^b$.
Choose a coefficient field
$\kappa_i\hookrightarrow R_i$, inducing the identity on the
residue field, by the Cohen structure theorem
\cite[Tag~032A]{Stacks}. We use the same symbol $\kappa_i$ for the
residue field and its chosen copy in $R_i$. This choice makes $R_i/\m_i^b$ a finite-dimensional
$\kappa_i$-algebra, and hence makes $Y_{i,b}$ projective over $\kappa_i$.
An ideal containing
$\m_i^b\OO_{Y_i}$ is exactly the inverse image of an ideal on the
thickening, so no information about the ideals with this bound is
lost. The coefficient fields help us parameterize this collection
in the proof; the collection itself was defined without them. In positive
characteristic $p$, the Cohen-ring map factors through its residue field
because $p=0$ in $R_i$, giving the required coefficient field. The chosen copy of $\kappa_i$ need not
contain the given copy of $k$. We use it only to form the parameter space;
quasi-projectivity in the enveloping condition remains over $k$.

For an ideal $h_i\supset\m_i^b\OO_{Y_i}$, define its length-Hilbert
polynomial by the eventual values
\begin{equation}\label{eq:polynomial}
 P_{i,h_i}(t)=\sum_{q\geq0}(-1)^q
  \length_{R_i}H^q\bigl(Y_i,(\OO_{Y_i}/h_i)\otimes A_f^{\otimes t}\bigr).
\end{equation}
The modules have finite length since their support is over the closed
point. Their lengths equal their dimensions over the coefficient field:
every factor in a composition series is the residue field $\kappa_i$,
which has dimension one over the chosen coefficient field.
Consequently \eqref{eq:polynomial} is the Hilbert polynomial of the
corresponding closed subscheme of $Y_{i,b}$ with respect to $A_f$.
Its coefficients in the binomial basis $\binom{t}{j}$ are integers, its
degree is bounded by $\dim Y_{i,b}$, and its definition is independent
of the coefficient field.

We now parameterize ideals with finitely many prescribed polynomials for
each $i$; the intrinsic coefficient bound in Section~\ref{sec:selection}
will provide such a finite list. A sufficiently large
power $A_f^{e_i}$ is very ample on $Y_{i,b}$. The Hilbert-scheme case of Grothendieck's representability theorem
\cite[Theorem~3.2, p.~260, and p.~265]{Grothendieck} gives a projective
parameter space for each polynomial $P(e_i t)$. Substitution
$P(t)\mapsto P(e_i t)$ is injective, so these parameter spaces give
exactly the required ideals. Write $S_i$ for their finite disjoint union.
Taking inverse images of the universal subscheme ideals under
$\OO_{Y_i}\to\OO_{Y_{i,b}}$ identifies $S_i(\kappa_i)$ with the
corresponding actual ideals on $Y_i$.

Give each $S_i(\kappa_i)$ its subspace Zariski topology and put
\[
 \mathcal P=\prod_{i=1}^s S_i(\kappa_i)
\]
with the ordinary finite product topology. We take the product of the
sets of rational points with their individual topologies; this does not
require a common field for the $\kappa_i$. Points of $\mathcal P$ represent
tuples of ideals. In the family argument below we distinguish a
specified dense set $\mathcal A$ of enveloping tuples from its closure.
We call members of $\mathcal A$, and their coordinate projections,
\emph{actual parameters}.

\begin{lemma}\label{lem:parameter}
The topological space $\mathcal P$ is noetherian.
\end{lemma}

\begin{proof}
Each $S_i(\kappa_i)$ is a subspace of a noetherian space, hence
is noetherian. For two noetherian spaces $S,T$, use noetherian induction
on pairs of closed subsets $F\subset S$, $G\subset T$ to show that
every open $\Omega\subset F\times G$ is quasi-compact. From an open
cover of nonempty $\Omega$, choose a nonempty basic rectangle
$U\times V$ contained in one member. The remainder is contained in
\[
 ((F\setminus U)\times G)\ \cup\ (F\times(G\setminus V)).
\]
Each product has a strictly smaller closed factor, so induction gives
a finite subcover of its intersection with $\Omega$. Together with
the chosen member this covers $\Omega$. Taking $F=S$, $G=T$ proves
noetherianness of $S\times T$; induction gives the finite product.
\end{proof}

\subsection{A closed condition on parameters}\label{sec:bad-parameters}

Let $C\subset\mathcal P$ be an irreducible closed subspace and let
$\mathcal A\subset C$ be a dense set of enveloping tuples.
We use $C$ only as a topological space, through irreducibility and the
density of $\mathcal A$. Universal ideals and containment equations
live separately on the Hilbert schemes $S_i$; their closed loci pull
back to $C$ along the coordinate projections. No common residue field,
universal family over $C$, or generic point of $C$ is needed. In
particular, closure points need not be enveloping. The following sum
uses only the actual members of $\mathcal A$:
\begin{equation}\label{eq:family-sum}
 h_i=\sum_{a\in\mathcal A}h_{a,i},\qquad
 H=H((h_i)),\qquad T=\nBl_HY.
\end{equation}
These ideal sums are coherent and equal finite subsums. Set
$H_a=H((h_{a,i}))$. Then $H_a\subset H$, and
Lemma~\ref{lem:contraction} identifies their pullbacks to $Y_i$
with $h_{a,i}$ and $h_i$, respectively.

We keep the normalized blowup $T$ of the sum fixed and vary the
member $a$. The ideal $H\OO_T$ is invertible. Wherever
$H_a\OO_T=H\OO_T$, the member $H_a$ therefore also becomes
invertible, giving a map to its blowup. This map will let us bring
an enveloping neighborhood from that blowup to $T$. Our immediate
task is to find a member for which the equality holds near a
prescribed point of a fiber. To do so, we describe the parameters
for which equality fails.

\begin{lemma}\label{lem:bad-parameters}
For each closed point $\tau\in T_{\eta_i}$ there is a proper closed subset
$F_\tau\subset C$, depending only on the $i$th parameter, such that
for $a\in\mathcal A$ one has $a\notin F_\tau$ exactly when
$H_a\OO_T=H\OO_T$ on a neighborhood of $\tau$.
\end{lemma}

The proof has two parts. First, failure to generate the invertible
ideal at $\tau$ is equivalent to containment in one fixed ideal
$N_\tau$ on $Y_i$. Second, containment in this fixed ideal is a
closed condition in the Hilbert parameter space. Properness of the
bad locus will then follow from the definition of the sum $h_i$.

\begin{proof}
Set $T_i=T\times_X\Spec R_i=T\times_YY_i$. Its closed fiber is
$T_{\eta_i}$; write $\widetilde\tau$ for the point of $T_i$
corresponding to $\tau$, and let $y\in Y_i$ be its image.
The local homomorphism
$\OO_{T,\tau}\to\OO_{T_i,\widetilde\tau}$ is flat and local,
hence faithfully flat. We may therefore check equality of the two
ideals after this pullback. The point $\widetilde\tau$ is closed
in $T_i$, and $T_i\to Y_i$ is projective, so $y$ is closed as well.
The field $k(\widetilde\tau)=k(\tau)$ is finite over $\kappa_i$,
since $\tau$ is closed in a fiber of finite type over $\kappa_i$.
It is consequently finite over $k(y)$.

Write $A_\tau=\OO_{T_i,\widetilde\tau}$ and let
$\mathfrak n_\tau$ be its maximal ideal. By
Lemma~\ref{lem:contraction}, extending $h_{a,i}$ and $h_i$ to
$A_\tau$ gives $H_aA_\tau$ and $HA_\tau$, respectively.
Choose a generator $u$ of the invertible ideal $HA_\tau$.
There is then an ideal $J_a\subset A_\tau$ such that
\[
 HA_\tau=uA_\tau,\qquad H_aA_\tau=uJ_a.
\]
Multiplication by $u$ identifies $A_\tau$ with $HA_\tau$.
Thus $H_aA_\tau=HA_\tau$ precisely when $J_a=A_\tau$, or
equivalently when $J_a$ has nonzero image in
$A_\tau/\mathfrak n_\tau=k(\tau)$. This is the rank-one case
of Nakayama's lemma. In particular, the relevant test takes place
in the fiber of the invertible ideal,
\[
 L_\tau=HA_\tau/\mathfrak n_\tau HA_\tau,
\]
which is one-dimensional over $k(\tau)$. Even if $u$ vanishes
at $\widetilde\tau$ as a function, its class generates this fiber.

Consider the natural map
$\varphi\colon(h_i)_y\to L_\tau$. It is nonzero because
$(h_i)_y$ generates $HA_\tau$. The map kills
$\mathfrak m_y(h_i)_y$, so it defines a sheaf map
\[
 h_i\longrightarrow(\iota_y)_*L_\tau,\qquad
 \iota_y\colon\Spec k(y)\hookrightarrow Y_i.
\]
Let $N_\tau$ be its kernel. This is a coherent ideal, since
$k(\tau)/k(y)$ is finite and hence the target is a coherent
skyscraper sheaf. It is properly contained in $h_i$ and agrees
with $h_i$ away from $y$. The preceding fiber test now says
\[
 \begin{aligned}
 H_aA_\tau=HA_\tau
 &\quad\Longleftrightarrow\quad
       \varphi((h_{a,i})_y)\ne0\\
 &\quad\Longleftrightarrow\quad h_{a,i}\not\subset N_\tau.
 \end{aligned}
\]
Faithful flatness gives the same test at $\tau$ on $T$.
Equality at that stalk holds on a neighborhood of $\tau$,
because the quotient $H\OO_T/H_a\OO_T$ is coherent.

We next express $h_{a,i}\subset N_\tau$ by equations on $S_i$.
Every parameter ideal contains $\m_i^b\OO_{Y_i}$. If this
fixed ideal is not contained in $N_\tau$, no parameter can be
bad, and we take $F_\tau=\varnothing$. Otherwise $N_\tau$
defines a closed subscheme $Z_\tau\subset Y_{i,b}$. An ideal
$h_{p,i}$ represented by a rational parameter $p\in S_i(\kappa_i)$
is the inverse image of the ideal of a subscheme
$Z_p\subset Y_{i,b}$. Reversing ideal inclusion gives
\[
 h_{p,i}\subset N_\tau
 \quad\Longleftrightarrow\quad Z_\tau\subset Z_p.
\]
We have reduced the question to whether a varying subscheme
contains the fixed subscheme $Z_\tau$.

Here are equations for this condition. Take affine opens
$\Spec A^\circ\subset Y_{i,b}$ and $\Spec B^\circ\subset S_i$, and let $M_\tau$
be the coordinate ring of $Z_\tau\cap\Spec A^\circ$. Choose finitely
many generators $g_\nu$ of the universal ideal in
$A^\circ\otimes_{\kappa_i}B^\circ$, and restrict them to
$M_\tau\otimes_{\kappa_i}B^\circ$. Only finitely many vectors of the
$\kappa_i$-vector space $M_\tau$ occur in these expressions. Choosing
a basis $m_1,\ldots,m_N$ of their span, write the restrictions as
\[
 \overline g_\nu=\sum_{\ell=1}^N m_\ell\otimes b_{\nu\ell}.
\]
No finite-dimensionality of the whole ring $M_\tau$ is needed here.
The universal quotient is flat over $S_i$, so specialization of
the universal ideal gives the ideal of the fiber $Z_p$. Moreover,
the $m_\ell$ remain linearly independent after extension to any
parameter residue field. The specialized generators therefore
vanish on the base change of $Z_\tau\cap\Spec A^\circ$ exactly when
\[
 b_{\nu\ell}(p)=0\qquad\text{for every $\nu,\ell$}.
\]
These are closed equations on $\Spec B^\circ$. A finite affine cover
of $Y_{i,b}$ tests containment on the whole subscheme $Z_\tau$.
The resulting closed conditions agree on overlapping parameter
charts because each expresses that same containment.

Pull this closed locus back to $C$ along the $i$th coordinate
projection and call it $F_\tau$. The equivalences above give its
required meaning on $\mathcal A$. Finally, if every actual
parameter were in $F_\tau$, then every $h_{a,i}$ would be
contained in $N_\tau$, and hence
\[
 h_i=\sum_{a\in\mathcal A}h_{a,i}\subset N_\tau.
\]
This contradicts the nonzero map $\varphi$. Thus at least one
actual parameter lies outside $F_\tau$, proving that it is proper.
\end{proof}

\subsection{An enveloping family sum}\label{sec:family-sum}

The preceding lemma treats one point at a time. Irreducibility
lets us choose a single member that works at any prescribed finite
set of points. We will transfer that member's quasi-projective
neighborhood to $T$, then use Lemma~\ref{lem:finite-neighborhoods}
to pass from finite sets to whole fibers.

\begin{proposition}\label{prop:pooling}
In \eqref{eq:family-sum}, the summed tuple $(h_i)$ is enveloping.
\end{proposition}

\begin{proof}
Let $\Sigma$ be any finite set of closed points of the fibers
$T_{\eta_i}$. For each $\tau\in\Sigma$,
Lemma~\ref{lem:bad-parameters} gives a nonempty open subset
$C\setminus F_\tau$. Since $C$ is irreducible, their finite
intersection
\[
 G_\Sigma=\bigcap_{\tau\in\Sigma}(C\setminus F_\tau)
\]
is nonempty and open. Density of $\mathcal A$ supplies an actual
parameter $a\in\mathcal A\cap G_\Sigma$. This one member
satisfies $H_a\OO_T=H\OO_T$ near every point of $\Sigma$.

The full equality locus is the open subset
\[
 O_a=T\setminus\Supp(H\OO_T/H_a\OO_T).
\]
It contains $\Sigma$. It also contains the inverse image of $W$,
which we identify with $W$, because both ideals are unit there.
On $O_a$ the ideal $H_a\OO_T$ is invertible. The universal
property of the blowup therefore gives a $Y$-morphism
$O_a\to\Bl_{H_a}Y$. The source is an open in the normal
integral variety $T$, and this morphism is dominant, since it is
the identity over the common generic point. It consequently
factors through the normalization, giving
\[
 g_a\colon O_a\longrightarrow T_a=\nBl_{H_a}Y.
\]

Because $a$ is enveloping, there is a $k$-quasi-projective open
$Q_a\subset T_a$ containing $W$ and every fiber $(T_a)_{\eta_i}$.
Set $U_a=g_a^{-1}(Q_a)$. This open contains $W$ and $\Sigma$:
the map is over $Y$, and hence over $X$, so a chosen point over
$\eta_i$ maps into the corresponding fiber of $T_a$.

We check that $U_a$ is quasi-projective over $k$. The morphism
$O_a\to Y$ is quasi-projective, since $O_a$ is open in the
projective $Y$-variety $T$. Factor $g_a$ as its graph followed
by projection:
\[
 O_a\xrightarrow{\ \Gamma_{g_a}\ }O_a\times_YT_a
       \xrightarrow{\ \operatorname{pr}_2\ }T_a.
\]
The graph is a closed immersion because $T_a\to Y$ is separated.
The projection is the base change of $O_a\to Y$, so it is
quasi-projective. Thus $g_a$ is quasi-projective. Restricting
over $Q_a$ and composing with $Q_a\to\Spec k$ shows that
$U_a$ is quasi-projective over $k$.

We have proved that every finite set $\Sigma$ lies, together
with $W$, in a quasi-projective open of $T$; the member $a$
and the open $U_a$ may depend on $\Sigma$. Lemma~\ref{lem:finite-neighborhoods}
applies because $T$ is normal and projective over $X$, and $W$
is quasi-projective. It gives a single quasi-projective open
containing $W$ and all the fibers $T_{\eta_i}$. The finiteness
of maximal quasi-projective opens is what makes this last passage
possible: if each maximal open containing $W$ missed a fiber
point, one missing closed point from each would form a finite
set that no such open could contain. The required neighborhood
therefore exists, so $(h_i)$ is enveloping.
\end{proof}

\subsection{A finite intrinsic collection}\label{sec:selection}

Return to the full poset of enveloping tuples for $f,W$, with no thickness
or polynomial bound imposed. For a maximal enveloping tuple $h$, let $b(h)$
be the least positive integer such that
\[
 \m_i^{b(h)}\OO_{Y_i}\subset h_i\qquad(1\leq i\leq s).
\]
We call $b(h)$ the \emph{thickness} of $h$.
At fixed thickness, the closed subschemes of $Y_{i,b}$ may still
have infinitely many Hilbert polynomials. The coefficient bound
below leaves only finitely many of them.
Write the polynomials \eqref{eq:polynomial} in the binomial basis and
let $c(h)$ be the largest absolute value of their coefficients, for all
$i$ and all degrees. Set
\[
 b_W=\min_h b(h),\qquad
 c_W=\min_{h:\,b(h)=b_W}c(h),
\]
where both minima range over maximal enveloping tuples. Define
\begin{equation}\label{eq:selected}
 \MM_W=\{h:\ h\text{ is maximal enveloping},\ b(h)=b_W,\ c(h)=c_W\}.
\end{equation}
Lemma~\ref{lem:existence} makes these minima well defined and the set
nonempty. Everything in \eqref{eq:selected} is intrinsic to $f,W$.

\begin{proposition}\label{prop:finite}
The set $\MM_W$ is finite.
\end{proposition}

\begin{proof}
At fixed thickness the dimensions of the thickenings $Y_{i,b_W}$ are
bounded. The coefficient bound therefore leaves finitely many possible
Hilbert polynomials. Section~\ref{sec:parameters}, with $b=b_W$ and this finite list
of polynomials, places $\MM_W$ in $\mathcal P$, which is noetherian
by Lemma~\ref{lem:parameter}.
Its closure has finitely many irreducible components $C_1,\ldots,C_N$.
Finite type alone would still allow infinitely many selected points
on a component; we show that maximality leaves at most one.

For every $j$, the subset $\mathcal A_j=\MM_W\cap C_j$ is dense
in $C_j$. Indeed, the complement in $C_j$ of the other components
is dense and open in the total closure, and every neighborhood of
each of its points meets $\MM_W$. Proposition~\ref{prop:pooling}
shows that the componentwise sum of the tuples in $\mathcal A_j$
is enveloping. It contains each of them, so full maximality makes each
member equal to that sum. There is at most one distinct tuple in
$\mathcal A_j$, proving finiteness.
\end{proof}

This proof explains why we took maxima before imposing the bounds.
The sum of a family may leave the numerical box, but it remains
in the full poset of enveloping tuples. Maximality in that poset
still applies. Had we required maximality only within the bounded
parameter family, this last inference would not follow.

\begin{remark}\label{rem:scheme-comparison}
One obtains exactly the same selected ideals by using
$\OO_{X,\eta_i}$ in place of its completion. Every ideal under consideration
contains a power of the maximal ideal, so both collections are inverse
images of ideals on the identical finite thickenings $Y_{i,b}$. This
identification preserves inclusions, sums, unit conditions and canonical
contractions, and therefore preserves the enveloping predicate and all
its maximal elements. The least thickness agrees. Flat cohomology base
change identifies the finite-length modules over the common Artinian
quotient, so their lengths and Hilbert polynomials agree as well.
Consequently the least pair $(b_W,c_W)$, the entire set $\MM_W$ and its
contracted ideals are unchanged. With the same choice of normalized blowup or principalization at each
round, the two rules agree throughout the tower.
\end{remark}

\section{The normal construction}\label{sec:construction}

The finite collections of Section~\ref{sec:selection} give us the
ideals for the next round. We use all of the maximal opens whose
boundary images have the smallest dimension, and take the product
of all their selected ideals. This makes the next round independent
of any choice among these opens or ideals.

For a retained tower $f\colon Y\to X$ with $Y$ not quasi-projective
and with the original mark $V\subset X$, define
\[
 \delta(f)=\min_{W\supset f^{-1}V}\dim f(Y\setminus W),
\]
where $W$ ranges over maximal quasi-projective opens of $Y$ containing
the displayed open. This is a nonnegative integer. Improving one open
attaining this minimum would suffice to lower $\delta$. To avoid
choosing that open or a witness for it, we use every minimizing open
and every selected tuple. Making their product invertible will allow
comparison with all the corresponding witnesses. For each minimizing
open, form the finite set \eqref{eq:selected}, and put
\begin{equation}\label{eq:round-ideal}
 J_f=\prod_{W:\,d_W=\delta(f)}\ \prod_{h\in\MM_W}H(h).
\end{equation}
This is a finite product of nonzero coherent ideals and is unit on
$f^{-1}V$. For the normal construction, one \emph{round} is the normalized blowup
\[
 p\colon Y'=\nBl_{J_f}Y\longrightarrow Y.
\]
We retain its full tautological Cartier divisor and append it to the
ordered tower. Since $Y$ is normal and $J_f$ is unit on $f^{-1}V$, the
new tower is still unchanged over $V$.

The next proposition applies more generally to any projective birational
map $p\colon Y'\to Y$ from a normal variety on which $J_f$ becomes an
invertible ideal, provided $p$ is unchanged over $f^{-1}V$. This will
allow us to use the same argument for smooth principalization.

\begin{proposition}\label{prop:drop}
Let $p\colon Y'\to Y$ be as in the preceding paragraph. Either $Y'$
is quasi-projective or $\delta(f\circ p)<\delta(f)$.
\end{proposition}

\begin{proof}
Fix one factor $H(h)$ of \eqref{eq:round-ideal}, associated to an open
$W$ with $d_W=\delta(f)$. Its total transform on $Y'$ is invertible.
To see this, work in a local domain where the product of the transformed
ideals is principal. Each factor has an inverse fractional ideal,
obtained by dividing the product of the other factors by that generator.
Thus it is invertible. Since $Y'$ is normal, the induced map to
$\Bl_{H(h)}Y$ factors through $T_h=\nBl_{H(h)}Y$.

The morphism $Y'\to T_h$ is projective. Its graph is
closed in $Y'\times_YT_h$, which is projective over $T_h$.
Let $Q_h\subset T_h$ be a witnessing open. Its inverse image in $Y'$ is
quasi-projective, contains $p^{-1}W$, and contains every fiber of
$Y'\to X$ over the top-dimensional generic points $\eta_i$ attached to $W$. Enlarge this inverse image to a maximal quasi-projective open
$W'\subset Y'$.

The closed image
\[
 (f\circ p)(Y'\setminus W')\subset f(Y\setminus W)
\]
misses every top-dimensional generic point of the right side. Its
dimension is therefore strictly less than $d_W$. If $Y'$ is
not quasi-projective, $W'$ is one of the opens in the definition of
$\delta(f\circ p)$ and gives the asserted inequality. If the image is empty, then $W'=Y'$, so $Y'$ itself is quasi-projective.
\end{proof}

\begin{theorem}[Marked normal construction]\label{thm:marked}
Let $X$ be a normal variety over an arbitrary field $k$, and let
$V\subset X$ be a specified nonempty quasi-projective open. The rule
above gives a finite ordered tower of normalized blowups with normal
quasi-projective model $P(X,V)$, unchanged over $V$. It is empty when
$X$ is quasi-projective and uses at most $\dim X$ nonidentity rounds.
It is natural under semilinear isomorphisms of the pairs $(X,V)$, with
compatible lifts respecting identities and composition. If $X$ is
proper, $P(X,V)$ is projective over $k$.
\end{theorem}

\begin{proof}
Start at the identity of $X$ and stop whenever the current variety is
quasi-projective. Otherwise form the finite product
\eqref{eq:round-ideal} and take its normalized blowup. Every term in
this product is unit on the lifted mark, so the same original $V$
is preserved at every stage. Proposition~\ref{prop:drop} strictly
lowers the nonnegative integer $\delta$ unless quasi-projectivity has
been reached. Initially $\delta\leq\dim X-1$, since every boundary
image is a proper closed subset of $X$. A stage with $\delta=0$ is
followed by a quasi-projective model. Thus there are at most
$\dim X$ rounds. If $X$ is proper, the final model is both proper
and quasi-projective, hence projective.

We check that the rule selects the entire tower without choices.
A semilinear isomorphism of marked inputs transports their maximal
quasi-projective opens containing the mark. Inductively, an isomorphism
of retained ordered towers transports the Cartier divisors and the
least positive ample weights in \eqref{eq:relative-polarization}.
It therefore transports $A_f$, the boundary images and their generic
points, the finite infinitesimal neighborhoods, and their completions.
Contraction and the enveloping condition are preserved, so the full
maximal tuples and their intrinsic thickness and Hilbert coefficients
are preserved. Their minima and the finite selected sets are preserved
as well. Consequently the product $J_f$ is transported exactly.

Blowup and normalization give a compatible isomorphism of the next
round, and quasi-projectivity gives the same stopping condition.
The coefficient fields and Hilbert schemes in
Section~\ref{sec:families} served only to prove finiteness; none is
chosen by the rule. At every stage a lift is unique because it has
its prescribed value on the dense unchanged open and its target is
separated. Uniqueness gives the identity and composition laws.
\end{proof}

For a normal variety, take $V=U(X)$ to obtain the scheme case of
Theorem~\ref{thm:normal}. The marked formulation will matter for
algebraic spaces: after making a scheme model, we retain the original
core as the mark instead of replacing it by the core of that new model.

\section{Canonical schematization}\label{sec:schematization}

The normal construction so far starts with a scheme. An algebraic
space has \'etale scheme charts, but need not have a Zariski open cover
by schemes. Thus local calculations on charts do not by themselves
give a scheme model. To apply the normal construction, we first make
such a model while preserving every open already known to be a scheme. This preliminary construction has a
simpler selection argument: its defining property is closed under arbitrary
sums, so a single greatest ideal replaces the finite collection of
Section~\ref{sec:selection}.

Throughout this section, $S$ is a noetherian scheme and $X\to S$ is a
separated morphism of finite type from an algebraic space. Normality and
reducedness are not assumed. Let $W=X_{\mathrm{sch}}$ be the largest
schematic open of $X$. It is dense by \cite[Tag~06NH]{Stacks} and
quasi-compact because $X$ is noetherian. In this section $W$ denotes
this schematic locus, rather than a maximal quasi-projective open.

\begin{definition}\label{def:schematizing}
A coherent ideal $I\subseteq\OO_X$ is \emph{schematizing} if
$I|_W=\OO_W$ and the ordinary blowup $\Bl_I X$ is a scheme.
\end{definition}

For a scheme $X$, the unit ideal is schematizing. For a general algebraic
space, the condition on the source of the blowup is substantive. The
morphism $\Bl_I X\to X$ is always representable by schemes; this does
not say that $\Bl_I X$ itself is a scheme. Blowing up the unit ideal,
for instance, gives $X$ again.

\subsection{Existence preserving the entire schematic locus}

\begin{lemma}\label{lem:schematizing-existence}
The space $X$ has a schematizing ideal.
\end{lemma}

\begin{proof}
Choose proper compactifications $X\subset\overline X$ and
$W\subset\overline W$ over $S$, with $\overline X$ an algebraic space
and $\overline W$ a scheme. The first is supplied by
\cite[Theorem~1.2.1]{CLO}, and the second by scheme compactification
\cite[Tag~0F41]{Stacks}. We need both compactifications because the
final source must be a scheme, while the modification must be over $X$.

Take the schematic closure $\Gamma$ of the common-open graph of $W$ in
$\overline X\times_S\overline W$. Its two projections are proper and
are isomorphisms over the specified copies of $W$:
\[
\begin{tikzcd}[column sep=large]
 & \Gamma\arrow[dl]\arrow[dr] & \\
 \overline X && \overline W.
\end{tikzcd}
\]
Here is why the isomorphism assertion holds on the entire inverse image
of $W$. Over $W\subset\overline X$, the graph of
$W\to\overline W$ is already closed, since $\overline W\to S$ is
separated. Restriction of schematic closure to this open therefore adds
nothing to that graph. The same argument over $W\subset\overline W$
uses the graph of $W\to\overline X$. Thus neither projection has
additional points or scheme structure above $W$.

Both inverse images of $W$ in $\Gamma$ are the same diagonal copy of
$W$, so admissibility refers to the same open for both projections.
The refinement theorem below is stated for representable morphisms of
stacks, where representability means representability by algebraic
spaces. A morphism between algebraic spaces has this property; its
source need not yet be a scheme.

Apply Rydh's refinement theorem \cite[Theorem~3.1]{Rydh} to
$\Gamma\to\overline X$. It gives a sequence of $W$-admissible blowups
$\Gamma'\to\Gamma$ such that $\Gamma'\to\overline X$ is also a sequence of
$W$-admissible blowups. The composite $\Gamma'\to\overline W$ is still
proper and isomorphic over the whole inverse image of $W$. Apply the
same theorem to this composite. We obtain $W$-admissible blowups
$\Gamma''\to\Gamma'$ such that $\Gamma''\to\overline W$ is a sequence of
$W$-admissible blowups. Since $\overline W$ is a scheme, $\Gamma''$ is a
scheme.

The other composite $\Gamma''\to\overline X$ is a sequence of
$W$-admissible blowups as well. Restrict it to $X\subset\overline X$.
Its source $\Gamma''_X$ is an open subscheme of $\Gamma''$. A composite of admissible
blowups of a noetherian algebraic space is one admissible blowup
\cite[Tag~086C]{Stacks}. Thus $\Gamma''_X=\Bl_I X$ for a coherent ideal
$I$ unit on $W$, which is the required schematizing ideal.

The argument applies also when $X$ is nonreduced. Rydh's theorem uses a
proper modification of algebraic spaces that is an isomorphism over the
specified open; it does not require that open to be schematically dense.
The conclusion preserves the scheme structure over $W$; it makes no
claim to preserve nilpotents supported outside that open.
\end{proof}

\subsection{Sums and the greatest ideal}

\begin{lemma}\label{lem:schematizing-sums}
The sum of any nonempty collection of schematizing ideals is schematizing.
\end{lemma}

\begin{proof}
We use the residual-ideal cover from
\cite[\S3, proof of Theorem~1.1 and Claim~3.3]{JiaMeng}.
For the present application a Zariski open cover by schemes suffices;
quasi-projectivity requires the additional finite-point argument of
Proposition~\ref{prop:pooling}.

Let $H=\sum_{I\in\mathcal S}I$, where every member of $\mathcal S$
is schematizing. Noetherianness and quasi-compactness imply that $H$ is
coherent and is a finite subsum. To see the finite-sum assertion, take
a finite affine \'etale cover of $X$, choose finitely many ideal members
whose sections generate the sum on each chart, and collect those finitely
many members. Equality with $H$ can be checked on the cover. Write
\[
 H=I_1+\cdots+I_m.
\]
This finite expression is used only in the proof; the sum itself is
independent of any such expression.

On $T=\Bl_H X$, the image ideal $L=H\OO_T$ is invertible. Since
$I_j\OO_T\subseteq L$, division by $L$ gives residual ideals
\[
 K_j=(I_j\OO_T)L^{-1}\subseteq\OO_T,
 \qquad K_1+\cdots+K_m=\OO_T.
\]
Their unit loci $O_j$ cover $T$. Indeed, at any local stalk at least one
$K_j$ is the unit ideal, since a sum of ideals contained in the maximal
ideal cannot be the whole ring. On $O_j$ we have
$I_j\OO_T=L$, so $I_j\OO_T$ is invertible there. The universal property
of the blowup gives an $X$-morphism
\[
 O_j\longrightarrow T_j=\Bl_{I_j}X.
\]

Each $T_j$ is a scheme. To deduce the same for $O_j$, factor this map
through its graph in $O_j\times_X T_j$. The product is a scheme:
$O_j\to X$ is representable by schemes, being an open subspace of the
blowup $T\to X$, and we have base-changed it by the scheme $T_j$.
The graph is closed because $T_j\to X$ is separated. Hence $O_j$ is
a scheme. The Zariski open cover by the $O_j$ now shows that $T$ is a
scheme. Finally, every $I_j$ is unit on $W$, so $H$ is unit there too.
\end{proof}

\begin{proof}[Proof of Theorem~\ref{thm:schematization}]
Let $\mathcal S$ be the nonempty collection of all schematizing ideals,
and put
\[
 H_X=\sum_{I\in\mathcal S}I.
\]
The two lemmas show that $H_X$ is itself schematizing. By construction,
it contains every member of $\mathcal S$, so it is the greatest one.
Its blowup is projective over $X$ and unchanged over $W$.

An isomorphism of the underlying algebraic spaces, without any
compatibility requirement on their chosen base presentations, transports
their largest schematic opens and their whole collections of
schematizing ideals. It therefore
transports $H_X$ itself. The resulting isomorphism of Rees algebras gives
a distinguished lift to the blowups; these lifts respect identities and
composition. This description depends only on $X$, not on its chosen
presentation over $S$. When $X$ is a scheme, $W=X$ and $H_X=\OO_X$,
so the construction is the identity.
\end{proof}

We call $\Bl_{H_X}X\to X$ the \emph{canonical schematization} and its
source the canonical scheme model. Greatestness here concerns the ideal
$H_X$ under inclusion. It does not assert a least model or a categorical
universal property. Over a general base $S$, the scheme model is projective
over $X$; it need not be projective over $S$.

The difference from Proposition~\ref{prop:pooling} is useful. Being a
scheme can be checked on the open cover by the $O_j$. Quasi-projectivity
cannot be inferred merely by covering with quasi-projective opens. The
main selection argument therefore needs one neighborhood containing all
the prescribed fibers, obtained by finite-point avoidance and the
finiteness of maximal quasi-projective opens.

\subsection{The normal model of an algebraic space}

\begin{proof}[Proof of Theorem~\ref{thm:normal}]
For scheme input, apply Theorem~\ref{thm:marked} with $V=U(X)$.
Now let $X$ be a normal integral algebraic space of finite type over $k$.
By Section~\ref{sec:prelim}, $U(X)$ is a nonempty quasi-projective open
contained in $X_{\mathrm{sch}}$. Form the canonical schematization and
normalize its source:
\[
 Z=\nBl_{H_X}X\longrightarrow X.
\]
This is a projective morphism from a normal scheme. Normalization is
finite here because the ordinary blowup is of finite type over a field.
Over $X_{\mathrm{sch}}$ the ordinary blowup is the identity, and that
open is already normal. The normalization is therefore also the identity
there. In particular, $U(X)$ identifies with a specified quasi-projective
open of $Z$.

Apply Theorem~\ref{thm:marked} to this marked pair $(Z,U(X))$.
Its original mark is preserved throughout. Composing with $Z\to X$
gives the required tower, with normal quasi-projective model $P(X)$.
The extra initial normalized blowup adds at most one round, and
$\dim Z=\dim X$. This bound reflects the preliminary schematization;
no optimality is asserted. If $X$ is proper, its model is proper and
quasi-projective, hence projective. If $X$ is quasi-projective, the
schematization and the marked tower are both the identity.

Both constructions transport their actual ideals under isomorphisms.
An isomorphism of fields transports the same definitions, so their
composition is natural under semilinear isomorphisms. Compatibility
through the ordered tower follows from the corresponding compatibility
in Theorems~\ref{thm:schematization} and~\ref{thm:marked}.
We retain $U(X)$ as the mark on $Z$; no relation between $U(Z)$ and
$U(X)$ is needed.
\end{proof}

\section{Smooth projective models in characteristic zero}\label{sec:smooth}

Assume now that the ground field has characteristic zero and the input
is smooth and proper. The same finite selection of enveloping tuples
can be followed by principalization of the selected product. This gives
smooth centers and a smooth projective model. We specify the retained
smooth tower first and prove the scheme case using
Sections~\ref{sec:prelim}--\ref{sec:construction}. For algebraic spaces
we extend this same selection directly, using \'etale descent of
principalization and intrinsic completed neighborhoods on which the
local calculations take place on schemes. This route keeps the smooth
construction independent of Theorem~\ref{thm:schematization} and uses
the same rule on schemes and spaces. The preliminary normal
schematization in Section~\ref{sec:schematization} is not used here.

\subsection{The fixed principalization procedure}\label{sec:principalizer}

For a nonzero ideal $I$ on a smooth variety $Y$, fix the procedure
obtained from the complete canonical resolution of the marked ideal
$(Y,I,\varnothing,1)$ over arbitrary fields of characteristic zero
\cite[Theorems~1.0.1, 2.4.1 and~2.4.2, \S\S5.3, 5.5 and~5.7]{Wlodarczyk2005}.
The chosen procedure is part of the construction; all canonicity
statements for $P_{\mathrm{sm}}$ are relative to this fixed choice.
It makes
a nonzero ideal invertible by an ordered finite sequence of smooth
blowups outside its unit locus. Indeed, each center lies in the support
of the current marked ideal, and that support lies over $V(I)$
\cite[Definitions~2.1.2--2.1.3 and \S5.5]{Wlodarczyk2005}.
Compatibility with semilinear
isomorphisms follows by transporting derivations: a semilinear
$\theta\colon A\to A'$ sends $D$ to $\theta D\theta^{-1}$, so
it transports derivative ideals. The remaining invariants and
centers are transported as well
\cite[Lemma~5.3.1 and Proposition~5.3.2]{Wlodarczyk2005}.
For an isomorphism between differently named ground fields, first
transport the scalar structure along that isomorphism; the same
derivation argument then applies.

The relatively ample bundle $A_f$ used to measure Hilbert polynomials
will depend on the ordered exceptional divisors, so we specify which
steps the tower records. The principalizer produces a sequence with
auxiliary divisor data: a blowup along a Cartier divisor is an
isomorphism but still contributes a recorded divisor. We compute the complete sequence first, then omit
every step whose morphism is an isomorphism. Each remaining step keeps
its full center, including any Cartier components, and the exceptional
divisor defined by that center. The towers in Theorem~\ref{thm:smooth}
and all subsequent choices use this convention.

For the retained smooth tower, use the total transforms of these full
divisors in \eqref{eq:relative-polarization}. Their tautological negative
line bundles are relatively ample, so the same least-weight construction
gives $A_f$. This agrees with the ordered-divisor convention in the
smooth construction throughout.

\subsection{The smooth construction for varieties}

\begin{proof}[Proof of Theorem~\ref{thm:smooth} for varieties]
Start with $f=\mathrm{id}_X$ and retain the original mark $U(X)$.
If the current variety is projective, stop. Otherwise use
Sections~\ref{sec:germs}--\ref{sec:selection}, with the relative
polarization of the retained smooth tower, to form $J_f$. Principalize
this ideal by the fixed procedure and append the retained smooth
blowups. This is one round, which may contain more than one blowup.

Every center avoids the inverse image of $U(X)$ because $J_f$ is unit
there. On the resulting smooth variety the image of $J_f$ is invertible,
so Proposition~\ref{prop:drop} applies. The boundary-image dimension
strictly decreases unless the new variety is quasi-projective. Since
all the varieties are proper, this stopping condition is projectivity.
There are at most $\dim X$ rounds, each containing a finite sequence
of smooth blowups. This proves termination and preservation of the core.

The transport argument in Theorem~\ref{thm:marked} applies to the same
ideals, polynomials and numerical minima on the retained smooth tower.
The chosen principalizer then transports every ordered center of the
next round. Uniqueness on the dense unchanged open gives compatibility
with identities and composition. The resulting smooth projective model
is denoted $P_{\mathrm{sm}}(X)$.
\end{proof}

The normal and smooth constructions use the same selection mechanism,
but retain different towers after a round. Their subsequent relative
polarizations and selected ideals need not agree. The smooth conclusion
comes from the principalization procedure, not from smoothness of
normalized blowups.

\subsection{Principalization on algebraic spaces}

For a smooth proper integral algebraic space, the maximal-open finiteness
and the finite-neighborhood lemma were established in
Section~\ref{sec:prelim}. Rydh's open-preserving Chow lemma applies to
algebraic spaces as well \cite[Theorem~3.5]{Rydh}. Combining its blowup
sequence by \cite[Tag~086C]{Stacks} gives the ideal form of
Lemma~\ref{lem:chow-ideal} with a projective scheme source. Thus the
witness needed for enveloping tuples is available in this category.

We first explain the compatibility needed for descent. Smooth
functoriality of the resolution sequence allows inserted isomorphism
steps \cite[Definition~2.1.5 and Theorem~2.4.2]{Wlodarczyk2005}.
For a surjective \'etale map we can align the actual centers, including
the auxiliary divisor data: the invariant values pull back unchanged,
and surjectivity preserves their global maximum. The maximum-center
construction therefore selects the inverse image of the same center.
The ordered exceptional divisors and their labels pull back as well.
Induction gives compatibility of the complete ordered sequences
\cite[Theorem~4.0.1(2),(4), Lemma~4.0.2 and \S\S5.1, 5.7]{Wlodarczyk2005}.
Surjectivity is useful
here even when the chart has several connected components: no invariant
value occurring on the target disappears on the chart.

Now choose an \'etale surjection $Z\to Y$ with $Z$ a finite disjoint
union of affine schemes. Since $Y$ is separated, $Z\times_Y Z$ is
a scheme, and its two projections are surjective \'etale.
Principalize on $Z$. The two pullbacks of the first center ideal to
$Z\times_Y Z$ agree, so descent of quasi-coherent ideals gives a
center on $Y$. Smoothness of the center and avoidance of the unit
locus are \'etale local. Blowup commutes with flat pullback, so we
repeat this argument on each subsequent stage. The same compatibility
on the triple fiber product gives the cocycle condition. The sequence
is finite because the sequence on $Z$ is finite, and common refinements
give independence of the chart.

Descend the complete sequence with its auxiliary divisor data before
omitting entire isomorphism steps as in Section~\ref{sec:principalizer}.
On the retained tower, keep the full Cartier divisor defined by each
retained center ideal, take its total transforms $E_j$, and use the
least weights in \eqref{eq:relative-polarization} to define $A_f$.
For existence of the weights, pull the retained tower back to a
surjective \'etale chart of $X$ that is a finite disjoint union of affine
schemes. The pulled-back tower consists of schemes, so
\cite[Tag~0C4K]{Stacks} applies there; relative ampleness is \'etale local
on the base. This gives exactly the same procedure and bundle on varieties.
The principalizer has this \'etale compatibility; the projectivization
rule is asserted to commute only with semilinear isomorphisms.

\subsection{Intrinsic completed neighborhoods}

Let $D$ be the reduced closure of a generic boundary point $\eta$. Although
$X$ need not have a schematic neighborhood of $\eta$, the finite thickenings
of $D$ do near $\eta$. Their local Artinian rings give the completion needed for the
construction. The following lemma supplies the local-ring facts used
in Section~\ref{sec:germs}. Its finite-level identity replaces
\eqref{eq:completion-thickening} in the contraction argument.

\begin{lemma}\label{lem:completed-point}
Let $X$ be a smooth integral algebraic space over $k$, let
$D\subsetneq X$ be a reduced integral closed subspace with generic
point $\eta$, and put $D_b=V(\mathcal I_D^b)$. The local rings
$\Lambda_b=\OO_{D_b,\eta}$ are well defined on schematic neighborhoods.
Their inverse limit $R$ is a complete noetherian regular local ring,
with residue field $\kappa=k(\eta)$ and $R/\m^b=\Lambda_b$.
There is a unique morphism $\alpha\colon\Spec R\to X$ inducing all
the finite-level maps $\Spec \Lambda_b\to X$. It is affine and flat, and
\begin{equation}\label{eq:finite-local-square}
 \Spec R\times_X\Spec \Lambda_b=\Spec \Lambda_b
\end{equation}
with the canonical projections. These constructions are transported
by semilinear isomorphisms.
\end{lemma}

\begin{proof}
The space $D$ has a dense schematic open containing $\eta$
\cite[Tag~06NH]{Stacks}. Its corresponding open in $D_b$ is a
nilpotent thickening of a scheme, hence is a scheme
\cite[Tag~05ZR]{Stacks}. The local ring $\Lambda_b$ is independent of
that open. It is noetherian with field reduction $\kappa$, and
therefore Artinian. The transition maps $\Lambda_{b+1}\to \Lambda_b$ are
surjective. Each $\Spec \Lambda_b\to X$ is a monomorphism, being a
localization in a schematic open of $D_b$ followed by its closed
immersion into $X$.

Since $X$ is separated, it is a decent algebraic space. We may
therefore choose an affine \'etale neighborhood $(U,u)\to(X,\eta)$
such that $u$ is the only point over $\eta$ and
$k(u)=\kappa$ \cite[Tags~03I7 and~0BBP]{Stacks}. Put
$B=\OO_{U,u}$, with maximal ideal $\mathfrak n$.
The \'etale pullback of $D$ is reduced, and its local ring at $u$
is the field $\kappa$, so $\mathcal I_DB=\mathfrak n$.

For each $b$, the scheme $U\times_X\Spec\Lambda_b$ is affine,
since $U$ is affine and $X$ is separated. It is \'etale of finite type
over the Artinian local ring $\Lambda_b$, and has the single point
$u$ with residue field $\kappa$. It is therefore
finite \'etale of rank one over $\Lambda_b$. Its local ring at $u$
is $B/\mathfrak n^b$: the map $\Spec\Lambda_b\to D_b$ is a
localization, and $\mathcal I_D^bB=\mathfrak n^b$.
A finite \'etale algebra of rank one is the base ring itself, so
the structural maps give compatible isomorphisms
\[
 \Lambda_b\xrightarrow{\sim}B/\mathfrak n^b.
\]
Passing to inverse limits identifies $R$ with $\widehat B$.
Since $U$ is smooth over $k$, the local ring $B$ is regular.
Its completion is a complete noetherian regular local ring, and
the displayed isomorphisms give $R/\m^b=\Lambda_b$.

Using $R=\widehat B$, define $\alpha$ as the composite
\[
 \Spec R\longrightarrow\Spec B\longrightarrow U\longrightarrow X.
\]
It induces the prescribed maps on all finite quotients and is flat,
since completion, localization and the \'etale map are flat.
It is affine because $X$ is separated: its base change by an affine
chart of $X$ is a closed subscheme of the product of two affine
$k$-schemes. To prove uniqueness, two maps inducing the given finite-level
maps have a closed equalizer in $\Spec R$. Its ideal is contained in
every $\m^b$, hence is zero. Thus the two maps agree.

The equality $\mathcal I_DB=\mathfrak n$ also gives
$\mathcal I_DR=\m$. Consequently
$\Spec R\times_XD_b=\Spec\Lambda_b$, with its specified map to $D_b$.
Since $\Spec\Lambda_b\to D_b$ is a monomorphism, this proves
\eqref{eq:finite-local-square}. The completed map $\alpha$ itself
need not be a monomorphism; the finite-level identity is what
Lemma~\ref{lem:contraction} uses.

Every point of $\Spec R$ maps to a generization of $\eta$, and only
the closed point lies over $\eta$. Indeed, points of $\Spec B$ are
generizations of $u$, the point $u$ is the only point of $U$ over
$\eta$, and $\mathfrak nR=\m$. This proves the cosupport assertion
used in Section~\ref{sec:germs} also for these completed neighborhoods.

The rings $\Lambda_b$, their inverse limit and the unique map $\alpha$
were defined from $X,D$. The \'etale chart only verifies their
properties. These intrinsic definitions give transport by semilinear
isomorphisms.
\end{proof}

\Needspace{4\baselineskip}
\subsection{The extended construction}

\begin{proof}[Proof of Theorem~\ref{thm:smooth} for algebraic spaces]
We explain separately the contraction, the parameter test, and the
transfer of neighborhoods. This shows where the scheme arguments apply
without changing the selection rule.

For a retained tower and a maximal open $W$, define $D_i,\eta_i$ as
in Section~\ref{sec:germs}. Lemma~\ref{lem:completed-point} gives
$R_i$ and the affine flat morphism $\Spec R_i\to X$. The pullbacks
$Y_i=Y\times_X\Spec R_i$ are projective schemes over $R_i$.
Their closed fibers and finite thickenings are exactly those in the
lemma. In particular, the cosupport of a germ unit on $W_i$ lies in the closed
fiber, so the germ contains a power of $\m_i$.

The affine flat maps $j_i\colon Y_i\to Y$ allow the contraction
\eqref{eq:contraction} on the algebraic space $Y$. Quasi-coherent
pushforward, flat pullback of kernels, and the finite-level identity
\eqref{eq:finite-local-square} give precisely the recovery argument of
Lemma~\ref{lem:contraction}. Coherence follows from noetherianness.
For the last assertion of that lemma, pull the two coherent difference
quotients to a surjective \'etale scheme chart $Z\to Y$. At any point
$z$ above one of the indicated fibers, their fibers over $k(z)$ vanish.
Nakayama's lemma gives zero stalks at $z$. The quotients also vanish
over the inverse image of $W$. Vanishing descends from the chart,
and coherence supplies the required open neighborhood of $W$ and
those fibers on $Y$.

The enveloping condition uses normalized blowups of these contractions.
The open-preserving Chow witness described above proves that enveloping
tuples exist. Noetherianness then gives maximal tuples, exactly as for
schemes. The length-Hilbert polynomials and the parameter spaces in
Section~\ref{sec:families} are formed on the schemes $Y_i$ and
$Y_{i,b}$, so they need no descent to the algebraic space.

For clarity, consider the equality test in Lemma~\ref{lem:bad-parameters}.
Let $T$ be the normalized blowup of a family sum and let
$\tau\in T_{\eta_i}$ be closed in the fiber. The scheme
$T_i=T\times_X\Spec R_i$ has the corresponding point
$\widetilde\tau$ with the same residue field. Choose an \'etale scheme
neighborhood $E\to T$ and a point $e$ over $\tau$, and then a point
$v$ of $E_i=E\times_TT_i$ over both $e$ and $\widetilde\tau$.
The local maps
\[
 \OO_{E,e}\longrightarrow\OO_{E_i,v}
 \longleftarrow\OO_{T_i,\widetilde\tau}
\]
are faithfully flat: the first is flat and local, and the second is
\'etale and local. Equality of coherent image ideals can therefore be
tested on either side. In particular, equality on $T$ near $\tau$ is
equivalent to equality on $T_i$ at $\widetilde\tau$; once it holds,
the complement of the support of the difference quotient is the
required open neighborhood.

Now the rank-one generation test takes place entirely in
$\OO_{T_i,\widetilde\tau}$. The image point $y\in Y_i$ is closed,
and $k(\tau)/k(y)$ is finite, as in Lemma~\ref{lem:bad-parameters}.
Thus the kernel $N_\tau$ is a coherent ideal on the scheme $Y_i$.
The fixed-ideal containment equations on the Hilbert scheme are
unchanged. We obtain the same proper closed bad-parameter loci and
can avoid any finite collection of them simultaneously. The auxiliary
\'etale neighborhoods verify the equality test; they do not enter the
parameter spaces or the selected ideals.

Finally, the graph argument of Proposition~\ref{prop:pooling} applies
to these algebraic spaces. The morphism $O_a\to Y$ is representable by
schemes and quasi-projective, since $O_a$ is an open of the projective
$Y$-space $T$. The graph of $g_a\colon O_a\to T_a$ is closed, so
$g_a$ is likewise representable by schemes and quasi-projective.
Over the witnessing quasi-projective open $Q_a\subset T_a$, its
source $U_a$ is consequently a quasi-projective scheme. Lemma~\ref{lem:finite-neighborhoods}
was already proved for normal algebraic spaces; it changes these
finite-point neighborhoods into one neighborhood of the whole fibers.
This proves the family-sum result and hence finite intrinsic selection.

Principalize the resulting product $J_f$ by the descended procedure.
Invertibility of every factor can be checked in \'etale local domains.
The graph comparison in Proposition~\ref{prop:drop} then gives a
quasi-projective open above each selected witnessing open. Properness
of the composite map makes its complementary image closed, and that
image misses every selected top-dimensional generic point. Thus the
same boundary-image dimension drops. The scheme-case proof gives
termination, preservation of the original core and semilinear
naturality. The resulting model is proper and quasi-projective, hence
a projective scheme, and is smooth and integral. Projective input
receives the empty tower.

On varieties, Lemma~\ref{lem:completed-point} recovers
$R_i=\widehat{\OO}_{X,\eta_i}$. The contractions, maximal enveloping
posets and polynomial bounds agree by
Remark~\ref{rem:scheme-comparison}; the full principalizer and divisor
conventions agree as well. Thus the two rules coincide at every round.
\end{proof}

\subsection{Moishezon manifolds}

The Artin--Moishezon correspondence supplies algebraic spaces even for
compact complex manifolds that are not algebraizations of varieties.
The construction below uses the fixed principalizer and divisor
convention of Section~\ref{sec:principalizer}.

\begin{corollary}\label{cor:moishezon}
Let $M$ be a compact connected Moishezon manifold. There is a finite
sequence of blowups in smooth centers
\[
 M_r\longrightarrow\cdots\longrightarrow M_0=M
\]
with $M_r$ projective, natural under biholomorphisms. If $X$ is the
proper algebraic space corresponding to $M$, the sequence is an
isomorphism over $U(X)^{\mathrm{an}}$. For the smooth projective
variety $P$ with $P^{\mathrm{an}}=M_r$, lifting gives an injection
\[
 \Aut_{\mathrm{hol}}(M)\lhook\joinrel\longrightarrow\Aut_{\mathbb C}(P).
\]
The modification is equivariant for this action, and $\dim P=\dim M$.
\end{corollary}

\begin{proof}
The classical Artin--Moishezon correspondence identifies proper
complex algebraic spaces with Moishezon spaces; for the precise
categorical statement and references to Artin's original work, see
\cite[\S1.1, p.~1]{Conrad}. Let $X$ correspond
to $M$. It is smooth: an \'etale chart analytifies locally
biholomorphically, and algebraic and analytic local rings at closed
points have the same completion. Connectedness makes $X$ integral,
since the irreducible components of a smooth algebraic space are
disjoint and open.

Apply Theorem~\ref{thm:smooth}. Analytification commutes with these
blowups: in \'etale-local coordinates with center $(z_1,\ldots,z_c)$,
both blowups have incidence equations $z_it_j=z_jt_i$, and the local
identifications glue. The final manifold is the analytification of the smooth projective
model $P=X_r$. Full faithfulness algebraizes
biholomorphisms and their inverses, so Corollary~\ref{cor:automorphisms}
gives the injection. The same equivalence and the tower's naturality
remove the algebraization choice and give naturality under
biholomorphisms.
\end{proof}

\begin{remark}[The equivariant-flattening assertion]\label{rem:kim-gap}
The proof of Theorem~1.1 in \cite[p.~564]{Kim} takes $X$ compact in
Fujiki's class $\mathcal C$ and chooses a surjective morphism
$\rho\colon X_1\to X$ from a compact K\"ahler space. It invokes
flattening to obtain a square
\[
\begin{tikzcd}
Y_1\arrow[r,"\sigma_1"]\arrow[d,"\pi"'] & X_1\arrow[d,"\rho"]\\
Y\arrow[r,"\sigma"'] & X
\end{tikzcd}
\]
in which the horizontal maps are projective modifications, $Y_1$ and $Y$
are smooth, $\pi$ is flat, and every automorphism of $X$ lifts through
$\sigma$. The asserted lifting condition is explicit:
``$\Aut(X)$ (resp.\ $\Aut(X_1)$) lifts uniquely and holomorphically
to $Y$ (resp.\ $Y_1$) by $\sigma$ (resp.\ $\sigma_1$).''
These conditions cannot be imposed for an arbitrary chosen $\rho$.

Take $\rho\colon\Bl_p\PP^2\to\PP^2$. This is a surjective morphism
from a smooth projective, hence compact K\"ahler, surface to a compact
space in class $\mathcal C$, so it satisfies the stated hypotheses on
the initial map. If every element of
$\operatorname{PGL}_3(\mathbb C)$ lifts through $\sigma$, its maximal
isomorphism open is nonempty and invariant under this transitive group.
It is therefore all of $\PP^2$, and $\sigma$ is an isomorphism.
After this identification, $\pi=\rho\circ\sigma_1$ is birational,
with zero-dimensional generic fiber. But $\sigma_1$ is surjective, so
$\pi^{-1}(p)$ maps onto the exceptional curve of $\rho$ and has positive
dimension. Flatness between the smooth surfaces would force all fibers
to have dimension zero, a contradiction.

This disproves the intermediate flattening assertion in the published
proof. It does not disprove its existential conclusion: the identity
of $\PP^2$ is an equivariant projective model. A specially compatible
choice of the initial K\"ahler map would require an additional argument;
equivariant resolution of the target alone does not provide that choice.
\end{remark}

\section{Comparison on Oda's threefold}\label{sec:oda}

Over $\mathbb C$, Oda's smooth complete nonprojective threefold
\cite[\S9, case~{$[7\text{-}5]$}]{Oda1978}
provides a comparison with the toric rule of
\cite[\S\S5--6]{BakhtaryToric}. Its fan has rays
\[
 \begin{gathered}
 v_1=e_1,\quad v_2=e_2,\quad v_3=e_3,\quad v_4=(-1,-1,-1),\\
 v_5=(-1,-1,0),\quad v_6=(0,-1,-1),\quad v_7=(-1,0,-1)
 \end{gathered}
\]
and maximal cones
\[
 123,\ 127,\ 136,\ 167,\ 235,\ 257,\ 356,\ 456,\ 457,\ 467,
\]
where $ijk$ denotes the cone generated by $v_i,v_j,v_k$.
Write $C_1=V(17)$, $C_2=V(25)$ and $C_3=V(36)$ for the invariant
curves corresponding to the indicated two-dimensional cones.
They are disjoint, since no maximal cone contains two of these walls,
and the rotation $(x,y,z)\mapsto(z,x,y)$ permutes them cyclically.
The three wall relations
\[
 v_2+v_6=v_1+v_7,\qquad
 v_3+v_7=v_2+v_5,\qquad
 v_1+v_5=v_3+v_6
\]
give
\begin{equation}\label{eq:oda-relation}
 L\cdot C_1+L\cdot C_2+L\cdot C_3=0
 \quad\text{for every line bundle }L\text{ on }X.
\end{equation}
In particular, $X$ is not projective. Blowing up $C_3$ makes it
projective. Here is the support-function certificate from
\cite[\S6, before Proposition~6.1]{BakhtaryToric}.
The blowup subdivides the cones $136$ and $356$ along their common
wall $36$, inserting $w=v_3+v_6=(0,-1,0)$. Define a function $h$ by giving the rays
$(v_1,\ldots,v_7,w)$ the values
\[
 (0,0,-8,-28,-18,-24,-22,-15),
\]
and extend linearly on each cone of the subdivided fan. If $\ell_\sigma$
is the resulting linear form on a maximal cone $\sigma$, these values
satisfy $\ell_\sigma(v)>h(v)$ for every ray $v$ outside $\sigma$;
the smallest difference is $1$. With this upper-convex convention,
the integral support function $h$ is strictly upper convex and certifies
projectivity. Rotation gives the same conclusion for $C_1,C_2$.

\subsection{The core and the selected ideal}\label{sec:oda-core}

Put $\Omega_i=X\setminus C_i$. These opens are quasi-projective because
they are unchanged in the projective curve blowups. To see that they
are all the maximal quasi-projective opens, let $W\subset X$ be any
quasi-projective open. Lemma~\ref{lem:chow-ideal} and principalization
give a smooth projective modification $g\colon P\to X$ that is an
isomorphism over $W$. Choose a very ample effective divisor $A$ on $P$
containing no exceptional prime divisor. The divisor $g_*A$ is effective and Cartier, since $X$ is smooth.
Its strict transform is $A$, so its pullback has the form
$g^*g_*A=A+F$ with $F$ effective and exceptional: the remaining
coefficients are the nonnegative vanishing orders of local equations
of $g_*A$ along exceptional divisors. If $C_i$ meets
$W$, its strict transform $\widetilde C_i$ is not contained in $F$, so
\[
 (g_*A)\cdot C_i=A\cdot\widetilde C_i+F\cdot\widetilde C_i>0.
\]
Equation~\eqref{eq:oda-relation} forces $W$ to miss some $C_i$ entirely.
Thus $W\subset \Omega_i$. The $\Omega_i$ are incomparable and hence are precisely
the maximal opens. It follows that
\[
 U(X)=X\setminus(C_1\cup C_2\cup C_3).
\]

For $f=\mathrm{id}_X$ and $\Omega_i$, the boundary image is the single
curve $C_i$. At its generic point $\eta_i$, the completed local ring
$R_i$ is regular of dimension two. Its maximal ideal $\m_i$ is
enveloping: its canonical contraction is $I_{C_i}$, whose blowup is
projective. The unit ideal is not enveloping, because a witnessing
open would strictly enlarge $\Omega_i$. Every proper ideal is contained
in $\m_i$, so this is the unique maximal enveloping ideal. Its
thickness and constant length-Hilbert polynomial are both $1$.
All three opens minimize the boundary-image dimension, giving
\begin{equation}\label{eq:oda-ideal}
 J_{\mathrm{id}_X}=I_{C_1}I_{C_2}I_{C_3}=I_{C_1\cup C_2\cup C_3}.
\end{equation}
Its ordinary blowup is the simultaneous blowup of the three disjoint
smooth curves. It is projective, being a further blowup of any one
of the projective curve blowups above. The fixed principalizer gives
a projective morphism to this blowup, so the rule stops after one round.
For the normal construction, the ordinary blowup is already smooth
and hence normal, so it is the canonical projective model in one
normalized round. For the smooth construction, this identifies the
selected ideal without specifying the smooth blowup sequence within
that round.

\begin{figure}[tb]
\centering
\begin{tikzpicture}[
  x=1.1cm,y=1.1cm,
  font=\small,
  line cap=round,line join=round,
  fan edge/.style={draw=black!75,line width=.55pt},
  selected wall/.style={draw=black,line width=1.6pt},
  new edge/.style={draw=black,line width=.75pt,dash pattern=on 2.1pt off 1.5pt},
  ray label/.style={fill=white,inner sep=.7pt,font=\footnotesize},
  panel title/.style={font=\small,align=center},
  face note/.style={font=\scriptsize,text=black!70}
]
\foreach \panel/\shift in {L/0,R/6.2}{
  \begin{scope}[xshift=\shift cm]
    \coordinate (\panel1) at (0,2.4);
    \coordinate (\panel2) at (-2.07846,-1.2);
    \coordinate (\panel3) at (2.07846,-1.2);
    \coordinate (\panel4) at (0,0);
    \coordinate (\panel5) at (0,-.76);
    \coordinate (\panel6) at (.65818,.38);
    \coordinate (\panel7) at (-.65818,.38);
    \fill[black!5] (\panel5)--(\panel6)--(\panel7)--cycle;
    \draw[fan edge] (\panel1)--(\panel2)--(\panel3)--cycle;
    \draw[fan edge] (\panel1)--(\panel6)--(\panel7)--(\panel2);
    \draw[fan edge] (\panel3)--(\panel5)--(\panel6);
    \draw[fan edge] (\panel5)--(\panel7);
    \draw[fan edge] (\panel4)--(\panel5)
                   (\panel4)--(\panel6)
                   (\panel4)--(\panel7);
    \draw[selected wall] (\panel1)--(\panel7)
                         (\panel2)--(\panel5)
                         (\panel3)--(\panel6);
    \node[face note] at (0,-1.62) {outer face $123$};
  \end{scope}
}

\coordinate (R8) at (-.32909+6.2/1.1,1.39);
\coordinate (R9) at (-1.03923+6.2/1.1,-.98);
\coordinate (R10) at (1.36832+6.2/1.1,-.41);
\draw[new edge] (R8)--(R2) (R8)--(R6)
                (R9)--(R3) (R9)--(R7)
                (R10)--(R1) (R10)--(R5);

\foreach \panel in {L,R}{
  \foreach \i in {1,...,7}{\fill (\panel\i) circle[radius=1.45pt];}
  \node[ray label,above=3pt] at (\panel1) {$v_1$};
  \node[ray label,below left=2pt] at (\panel2) {$v_2$};
  \node[ray label,below right=2pt] at (\panel3) {$v_3$};
  \node[ray label,above=3pt] at (\panel4) {$v_4$};
  \path (\panel5) ++(0,-.16) node[ray label] {$v_5$};
  \path (\panel6) ++(-.49,.17) node[ray label] {$v_6$};
  \path (\panel7) ++(.36,.17) node[ray label] {$v_7$};
}
\foreach \i in {8,9,10}{\draw[fill=white,line width=.8pt] (R\i) circle[radius=2.45pt];}
\path (R8) ++(.31,.12) node[ray label] {$v_8$};
\node[ray label,above left=3pt] at (R9) {$v_9$};
\path (R10) ++(-.24,-.27) node[ray label] {$v_{10}$};

\path (L1)--(L7) node[midway,ray label,right=3pt] {$C_1$};
\path (L2)--(L5) node[midway,ray label,above=3pt] {$C_2$};
\path (L3)--(L6) node[midway,ray label,right=3pt] {$C_3$};
\node[panel title] at (0,2.98) {(a) Oda's fan $\Sigma$};
\node[panel title] at (6.2/1.1,2.98) {(b) Fan of $\operatorname{Bl}_{C_1\cup C_2\cup C_3}X$};
\end{tikzpicture}
\caption{The fan of Oda's threefold and the fan of the ordinary blowup
of $I_{C_1}I_{C_2}I_{C_3}$. These are planar combinatorial drawings of
the intersections of the fans with a sphere, with $123$ the outer face.
Bold edges mark the walls $17$, $25$, $36$ and their subdivisions;
hollow vertices are the new rays $8$, $9$, $10$, and dashed segments
are the new walls. The shaded cones $456$, $457$, $467$ and the outer
face $123$ are exactly the maximal cones containing none of the three
selected walls. Their charts lie in $U(X)$ and are unchanged. The right-hand fan is the normal construction's
projective model, and is dominated by the smooth principalization round.}
\label{fig:oda}
\end{figure}
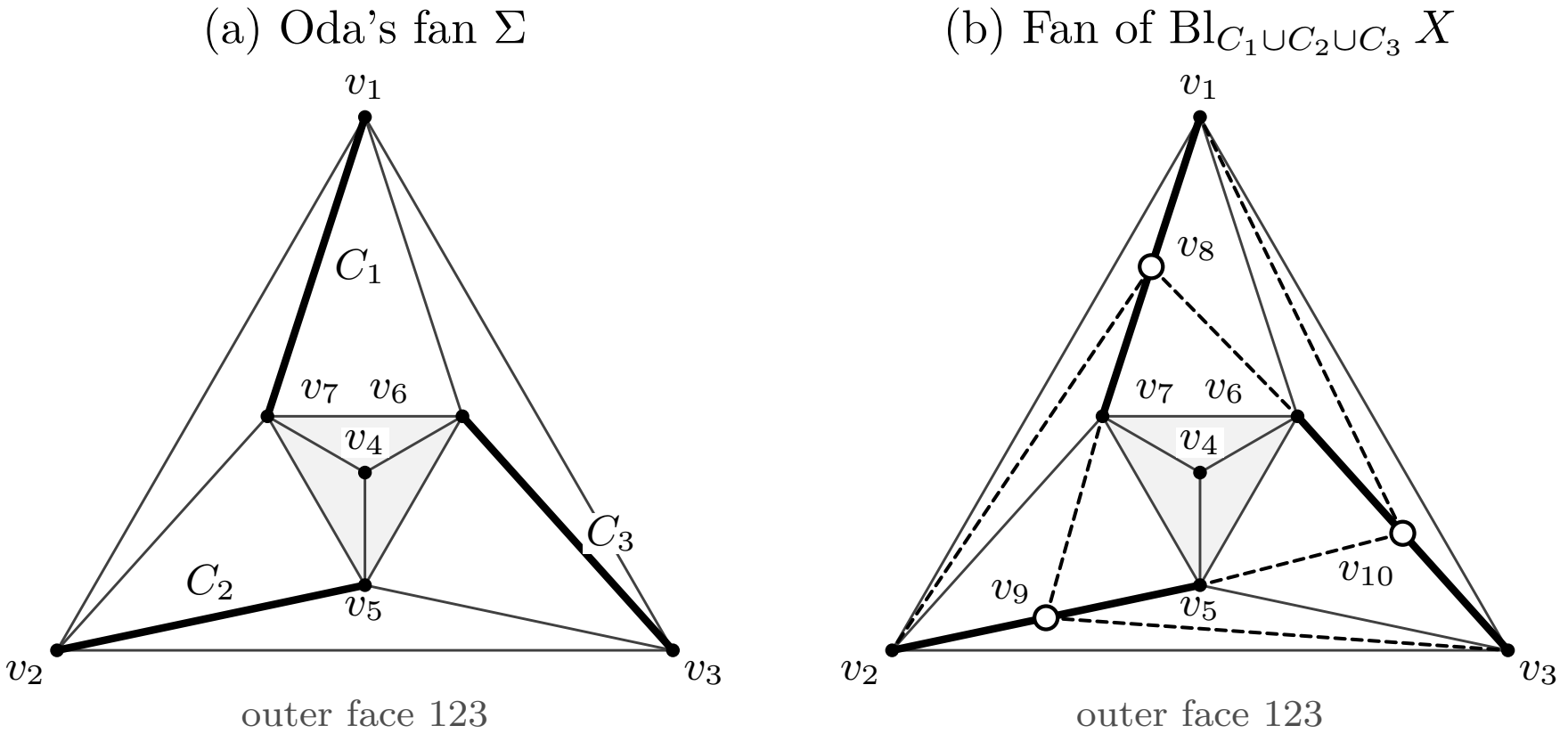

The fan operation is especially concrete. Blowing up the three invariant
curves inserts the rays
\[
 v_8=v_1+v_7=(0,0,-1),\qquad
 v_9=v_2+v_5=(-1,0,0),\qquad
 v_{10}=v_3+v_6=(0,-1,0).
\]
It subdivides the walls $17$, $25$ and $36$ and the maximal cones
containing them. Figure~\ref{fig:oda} displays the two triangulated
spheres using the same positions for the original rays.

\subsection{Comparison with the wall ideal}

The earlier canonical toric rule principalizes the product of
\[
 \mathcal J_{[m]}=\OO_X(-D_m^+)+\OO_X(-D_m^-)
\]
over all unoriented primitive wall characters $[m]$, where
$\operatorname{div}(\chi^m)=D_m^+-D_m^-$ is the decomposition into
zero and pole divisors \cite[\S5.1, equation~(9)]{BakhtaryToric}.
On the chart $U_{123}\cong\mathbb A^3$ with coordinates $x,y,z$,
the wall $67$ has normal $(1,1,-1)$ and its factor is $(xy,z)$.
This is nonprincipal at the origin. Since invertibility of a product
forces invertibility of every factor, as in Proposition~\ref{prop:drop},
principalizing the wall product must modify this chart.
In contrast, $U_{123}$ meets none of the $C_i$ and lies in $U(X)$;
our rule preserves the whole chart. Thus the two constructions give
different modifications over the same toric variety. The wall rule
makes every wall-character ideal invertible, whereas the present rule
preserves the original quasi-projective core.

For nonprojective $X$, Proposition~\ref{prop:drop} bounds the number
of rounds by $\delta(\mathrm{id}_X)+1\leq\dim X$. It gives no bound
on the number of nondivisorial connected components of the smooth centers. The threefolds $\Sigma(a,b)$ of
\cite[Theorem~3.2]{BakhtaryToric} require at least $\min(a,b)$ new rays
in every smooth projective refinement. Our rule gives such a refinement. Over $\mathbb C$, naturality
makes the round ideal $J_f$ invariant under every torus translation.
On an invariant affine chart, a torus-stable ideal is a sum of
character weight spaces, so it is generated by monomials. This gives
invariance under the algebraic torus itself. The fixed principalizer
preserves that invariance; its invariant centers give algebraic torus
actions on the blowups. A connected smooth invariant center is
irreducible and hence is an orbit closure. Each nondivisorial connected
component therefore adds one ray, so
at least $\min(a,b)$ such components occur in the whole tower.
This counts connected components of centers: a disconnected center
may add several rays in one blowup.

\subsection{Infinite action on cohomology}

The absence of an invariant-big-class hypothesis has content even for
smooth proper algebraic varieties. The following example makes this
visible using the same Oda threefold.

\begin{example}\label{ex:infinite-cohomology}
Let $Z$ be Oda's threefold above, let $E$ be a complex elliptic curve,
and put $X=Z\times E^2$. This is a smooth proper nonprojective fivefold:
if it were projective, the closed subvariety $Z\times\{a\}$ would
be projective. Its analytification is a compact Moishezon manifold.

Consider the automorphism
\[
 g(z,x,y)=(z,2x+y,x+y).
\]
The matrix
\[
 M=\begin{pmatrix}2&1\\1&1\end{pmatrix}
\]
has determinant $1$, so it defines an automorphism of $E^2$. Its
eigenvalues are $\lambda=(3+\sqrt5)/2$ and $\lambda^{-1}$.
Each occurs twice on $H^1(E^2,\mathbb R)$, since the coordinate action
is tensored with the two-dimensional space $H^1(E,\mathbb R)$.
The cohomology algebra of a torus is the exterior algebra on $H^1$.
Wedging two independent $\lambda$-eigenvectors therefore gives a
$\lambda^2$-eigenvector in $H^2(E^2,\mathbb R)$.

The projection $X\to E^2$ has a section $a\mapsto(z_0,a)$, so its
pullback on $H^2$ is injective and intertwines the indicated
automorphisms. Thus $g^*$ has eigenvalue $\lambda^2>1$ on
$H^2(X,\mathbb R)$, and its action on $H^2(X,\mathbb Q)$ has infinite
order.

There is no $g$-invariant big real $(1,1)$-class on $X^{\mathrm{an}}$.
Indeed, Jia--Meng's finite-index theorem
\cite[Corollary~1.3]{JiaMeng} would place some positive power of $g$
in the identity component of the automorphism group. That component
acts trivially on $H^2$, contradicting the eigenvalue $\lambda^{2m}>1$
for every $m>0$. In particular, no big class is invariant under the
full automorphism group. Theorem~\ref{thm:smooth} nevertheless supplies
a smooth projective model to which the full group lifts.
The example separates the hypotheses of the two results: its cyclic
action belongs to no big-class stabilizer. A projective model for this
particular cyclic action is already visible by taking the product of
the curve blowup of $Z$ with $E^2$, since $g$ acts trivially on $Z$.
The example is not intended to make that special lifting problem difficult.
\end{example}

\section*{Acknowledgements}

The author used OpenAI Codex and Anthropic Claude for proof development,
literature checks, and manuscript preparation and review. The author takes
full responsibility for the content.

\enlargethispage{2\baselineskip}

\end{document}